\documentclass[a4paper,11pt,reqno]{amsproc}

\usepackage{amsthm}
\usepackage{amsmath}
\usepackage{amssymb}

{ \theoremstyle{plain}
\newtheorem{theorem}{Theorem}
\newtheorem{lemma}{Lemma}
\newtheorem{corollary}{Corollary}
\newtheorem{proposition}{Proposition}
}

{ \theoremstyle{definition}
\newtheorem{definition}{{\bf Definition}}}

{ \theoremstyle{definition}
\newtheorem{example}{Example}
}

{ \theoremstyle{remark}
\newtheorem{remark}{Remark}
}

\def\U{{U}}
\def\J#1#2#3{\int\limits_{#3}|\nabla #2|^{#1}}
\newcommand{\lalfa}[2]{\lambda_{\p,#1}(#2)}
\newcommand{\grnd}[1]{\overline{\nabla}#1}
\newcommand{\mm}[1]{|#1|}
\newcommand{\treu}[1]{\bigtriangleup_{#1}}
\newcommand{\R}[1]{\mathbb{R}^{#1}}
\newcommand{\scal}[2]{{\mathcal h} #1, #2{\mathcal i}}

\newcommand{\M}{{\mathcal M}}
\newcommand{\D}{{\mathcal D}}
\newcommand{\OO}{{\mathcal O}}
\newcommand{\N}{{\mathcal N}}
\newcommand{\gd}[1]{\nabla #1}
\newcommand{\capac}{{\rm cap}_{\p}}
\newcommand{\et}{e^{\top}}
\newcommand{\en}{e^{\bot}}
\newcommand{\p}{\alpha}
\newcommand{\q}{V_2(\M)}
\def\Z{{\mathcal Z}}
\def\reduc{\lambda^\star}
\begin{document}

\title[Denjoy-Ahlfors Theorem]{Denjoy-Ahlfors Theorem for Harmonic Functions on
Riemannian Manifolds and External Structure of Minimal Surfaces }

\author{Vladimir Miklyukov and Vladimir Tkachev}
\thanks{This paper was supported by Russian Fundamental Researches
Fund, project 93-011-176}

%%%%%%%%%%%%%%%%%%%%%%%%%%%%%%%%%%%%%%%%%%%%%%%%%%%%%%%%%%%%%%%%%%%%%%%%%%%%

\begin{abstract}
We extend the well-known Denjoy-Ahlfors theorem on the
number of different asymptotic tracts of holomorphic functions to
subharmonic functions on arbitrary Riemannian manifolds. We obtain
some new versions of the Liouville theorem for $\p$-harmonic functions
without requiring the geodesic completeness requirement of a manifold.
Moreover, an upper estimate of the topological index of the height
function on a minimal surface in $\R{n}$ has been established and,
as a consequence, a new proof of Bernstein's theorem on entire solutions has been
derived. Other applications to  minimal surfaces are also discussed.
\end{abstract}
%\tableofcontents

\maketitle

\section*{Introduction}

The present paper is an attempt to develop efficient
methods for investigation of external geometry of higher-dimensional
minimal surfaces. The main difficulty of this theory is a lack of
representations  similar to the well known
Weierstrass representation for two-dimensional minimal surfaces. On the other hand, the coordinate functions of
a minimal immersion are harmonic  in the inner metric, and from
this point of view, the immersed minimal submanifold can be regarded as Riemannian manifolds $M$ equipped with a system of
harmonic functions $\{f_1(m),\ldots,f_n(m)\}$ satisfying certain structure
conditions like the following:
$$
\sum_{k=1}^{n}|\nabla f_k(m)|^2={\rm dim}\;M.
$$
Moreover, the minimality condition itself (i.e. the vanishing of the mean curvature) is a source for numerous
algebraic combinations between the coordinate functions of $M$ and
its Gaussian map which are subharmonic in the inner metric.

These properties together with suitable functional results about
subharmonic functions allows us to derive an information about the
geometrical structure of minimal immersions. We mention the substantial
papers of Yau \cite{Y}, Cheng and Yau \cite{C-Y}-\cite{CY} realizing
this way on complete higher-dimensional Riemannian manifolds with
restrictions on the Ricci curvature. The main tool was the maximum
principle for some classes of elliptic type PDE (see also \cite{Hilde},
\cite{Hilde1}). On the other hand, useful geometric applications of
the Liouville and Phragmen-Lindel\"off type theorems have been
obtained in \cite{Mik0}-\cite{M2}, \cite{MT1},\cite{MT2} by using
of methods  quasiconformal mapping theory. We notice that such
results should be regarded as corollaries of an appropriate form of
general Denjoy-Ahlfors type theorems on Riemannian manifolds.

The celebrated Denjoy-Ahlfors theorem states that
the order of growth of an entire holomorphic function $f(z)$ yields
an upper bound on the number $\rho$ of different asymptotic values
of $f(z)$. More precisely, let $w=f(z)$ be such a function. A family
of domains $\{\D(\tau)\}$ is called an \textit{asymptotic tract} of $w=f(z)$ if

a) every $\D(\tau)$ is an open component of the set
$$
\{z\in \mathbb{C}: |f(z)|>\tau>0\};
$$

b) for all $\tau_2>\tau_1>0$ there holds $\D(\tau_1)\supset\D(\tau_2)$
and $\cap_{\tau}\overline{\D}(\tau)=\emptyset$.

Two asymptotic tracts $\D'(\tau)$ and $\D''(\tau)$ are called different
if for a large enough $\tau>0$ we have
$\D'(\tau)\cap\D''(\tau)=\emptyset$.

Then the mentioned inequality is
$$
\limsup_{r\to+\infty}\frac{\log\log M(r)}{\log r}\geq
\frac{1}{2}\rho,
$$
where $M(r)$ is the maximum modulus of $f(z)$ on $|z|=r$.

It is well known \cite{HK} that this result gives an estimate
of the number of different {\it asymptotic tracts} of a
holomorphic function by its lower order. This terminology goes
back to \cite{MckLane} and seems to us more preferable by virtue of
its better suitability for geometrical applications.

In this paper we extend the Denjoy-Ahlfors theorem on Riemannian
manifolds. In addition, this geometrical approach enables us to
apply the developed technique  to minimal surface theory.
We describe briefly some key ideas of our method. Let $M$ be a noncompact
$p$-dimensional Riemannian manifold and let $(M,u)$ be a
$C^2$-surface given by an immersion $u(m):M\to\R{n}$. The
surface $\M\equiv(M,u)$ is called {\it minimal} if the mean
curvature vector is identically zero. If $\M$ is a minimal
surface then for an arbitrary vector $e\in\R{n}$ the
corresponding coordinate function $f(m)=\scal{e}{u(m)}$ is
harmonic with respect to the inner metric of $\M$ \cite{KN}.
Then asymptotic tracts of $f(m)$ can be identified with the
noncompact components which one can cut off from $u(M)$ by
hyperplanes $\scal{u}{e}={\rm const}$. On the other hand, there
exists a priori an upper bound on the growth of $f(m)$ in terms
of the growth of the distance function in $\R{n}$. Hence, a
Denjoy-Ahlfors type theorem applied to $f(m)$ on a minimal surface
$\M$ provides a relation between inner geometric
characteristics of $\M$ and the number of pieces which one can
cut off from this surface by hyperplanes.

The first three sections of the paper are devoted to several
generalizations of the classical Denjoy-Ahlfors theorem and their
corollaries on arbitrary, not necessary geodesically complete,
manifolds. In Section~\ref{sec3} we derive a first order
differential inequality on the Dirichlet integral of a subharmonic
function in terms of the $N$-means of fundamental frequency of its
level sets. This inequality depends on the topological structure of
asymptotic tracts of harmonic functions. We distinguish two different cases according to contain or not level sets of a harmonic
compact components. If harmonic functions contain compact level sets (on topologically nontrivial manifolds) then they generate
singular asymptotic tracts. The corresponding
theory is rather complicated because an non-zero additive term in the
main differential inequality. This case is closely linked with
minimal surfaces of tubular type. We treat singular tracts in Section~\ref{sec4}.

It was recently shown in \cite{T1} that there exists a finite or
infinite limit
$$
V_p(\M)=\frac{1}{\omega_{p}}\lim_{R\to
\infty}\frac{1}{\ln R} \int_{M_a(R)}\frac{dx}{|x-a|^p},
$$
where $M_a(R)$ is the part of the minimal surface $\M$ in the spherical shell $\{1<|x-a|<R\}$. $V_p(\M)$ is called
the {\it projective volume} of $\M$. It follows from \cite{T1} that
this quantity does not depend on the choice of $a$, and the number of
ends of $\M$ is estimated by $c(p,n) V_p(\M )$. Moreover, numerous connections between $V_p(\M)$ and the
integral-geometrical means of $\M$ have been established in \cite{T1}.

In the final section, we apply the developed methods to
the external structure of minimal surfaces with finite
projective volume.

\medskip
{\bf Theorem~\ref{th53}.}
\it
Let $\M $ be a two-dimensional properly immersed minimal surface
in $\R{n}$ of finite topological type, $e$ is a regular
direction and $x_1$ be the corresponding coordinate function. If
$V_2(\M)<+\infty$ then either $\M $ is a subset of a hyperplane $x_1
= {\rm const}$, or the number of critical points $\{ a_i \}$ of
$x_1(m)$ is finite. Moreover, the following estimate on the
topological index of $\M$ holds:
$$
\sum_{j} {\rm ind}(a_j) \leq \q - \chi (M).
\eqno{(\ref{pom})}
$$
Here $\chi (M)$ is the Eulerian characteristic of $M$.
\rm
\medskip

When $\M$ is homeomorphic to a sphere with a finite number of
points removed, Theorem~\ref{th53} was announced without proof in
\cite{M3}.

Let $\M$ be a two-dimensional simply connected surface. By $n(t
e^{i\theta})$ we denote the multiplicity of the orthogonal
projection of $\M$ onto a fixed twodimensional plane $V$ at the point
$z=te^{i\theta}$. The following result extends the celebrated Bernstein's theorem.

\medskip
{\bf Theorem~\ref{th54}.}\
\it
Let $\M$ be a two-dimensional minimal properly immersed plane in
$\R{3}$ properly projected onto $V$.
Then either $\M$ is a plane in $\R{3}$, or
$$
\liminf_{R\to\infty} \frac{1}{\ln R} \int_{1}^{R} \frac{dt}{t}
\int_{0}^{2\pi} n(t e^{i\theta}) d\theta \geq 8.
$$
\rm
\medskip

We establish the most of preliminary propositions in the form
which is applicable as well for $\p$-subharmonic functions. In particularly,
this allows us to apply these results to the so-called $\p$-minimal surfaces
(see the precise definition in Section~\ref{sec5}, cf.
\cite{TVG})

Finally, we agree upon the notations. We assume henceforth that $M$
is a $p$-dimensional noncompact orientiable Riemannian
$C^2$-manifold. By $\partial \OO$ we denote the boundary of a
(sub)manifold $\OO$. If $\partial M$ is a nonempty set we
assume that it is piecewise smooth. By $T_m(M)$, $\scal{X}{Y}$
and $\nabla$ we denote the tangent space to $M$ at $m\in M$, the
scalar product of two vectors $X,Y\in T_m(M)$ and Levi-Civita
connection on $M$ respectively.

Let $\M=(M,u)$ be a surface defined by a $C^2$-immersion
$u(m):M\to \R{n}$, $p<n$. By  $\nabla$ we denote the Riemannian
connection induced by the immersion $u(m)$. Let
$\overline{\nabla}$ denote the Euclidean covariant derivative in $\R{n}$
and $N_m(M)$ denote the normal space at $m$. Then, given $\xi\in
N_m(M)$, $X,Y\in T_m(M)$ we have the Gauss formulas
$$
(\overline{\nabla}_XY)^\top=\nabla_XY, \quad
(\overline{\nabla}_XY)^\bot =B(X,Y),
$$
and the Weingarten formulas
$$
(\overline{\nabla}_X \xi)^\top=-A^\xi (X), \quad
(\overline{\nabla}_X \xi)^\bot=\nabla_X\xi.
$$
Here the top symbols $\top$ and $\bot$ denote the orthogonal projections
on $T_m(M)$ and $N_m(M)$ respectively. Here,  $B$ and $A$ are
the second fundamental form and Weingarten map respectively. Then the
following relation holds
$$
\scal{A^\xi(X)}{Y}=\scal{B(X,Y)}{\xi}.
$$

In what follows we also skip the notation of the  Lebesgues measure in integrals.

\smallskip
{\it Acknowledgements}: \ The authors thank the referees for many
helpful suggestions that greatly improved the presentation of
this paper.

%%%%%%%%%%%%%%%%%%%%%%%%%%%%%%%%%%%%%%%%%%%%%%%%%%%%%%%%%%%%%%%%%%%%%%
\section{Subharmonic and exhausting functions}
%%%%%%%%%\setcounter{equation}{0}
%%%%%%%%%%%%%%%%%%%% 1st Section %%%%%%%%%%%%%%%%%%%%%%%%%%%%%%%%%%%%%

\subsection{} We recall that the divergence of a smooth vector field $X$ is
defined to be
$$
{\rm div} X= \sum_{i=1}^{p}\scal{\nabla_{E_i}X}{ E_i}
$$
where the last sum does not depend on choice of the orthonormal basis
$\{E_i\}$ of $T_m M$.

Let $\p >1$ be a fixed real parameter. We define the
differential operator
$$
\Delta_\p f = {\rm div} (|\nabla f|^{\p-2}\nabla f)=
|\nabla f|^{\p-3}(|\nabla f|\Delta f + (\p - 2)\scal{\nabla
f}{\nabla|\nabla f|})
%\label{ua}
$$
where $\Delta\equiv\Delta _{2}$ is the ordinary Laplace operator on
$M$.

We denote by ${\mathcal Z}(f)$ the set of all critical points of a
function $f(m)$.

\smallskip
\begin{definition} A smooth function $f(m)$ is called
$\p $-{\it subharmonic} if the inequality
\begin{equation}
\Delta _{\p} f(m) \geq 0
\label{ub}
\end{equation}
holds everywhere in $M\setminus {\mathcal Z}(f)$. If the boundary
$\partial M$ is nonempty we assume that everywhere on $\partial M$
holds the Neumann condition
\begin{equation}
\scal{\nabla f(m)}{\nu (m)} = 0
\label{Neimann}
\end{equation}
where $\nu (m)$ is the outward normal to the boundary $\partial M$.
\end{definition}

\begin{definition} If (\ref{ub}) turned into equality the
function $f(m)$ is called $\p$-{\it harmonic}.
\end{definition}

\subsection{} Let $h(m): M \rightarrow (0;h_0) $ be a continuous function
of $C^{\infty}(M \backslash \Sigma_0)$, where $\Sigma
_0$ is the null-level set of $h(m)$. We introduce for $t \in
(0;h_0)$ an $h$-ball as
$$ B_h(t)=\{ m\in M: h(m)<t\},
$$
and
an $h$-sphere as
$$
\Sigma _h (t)= \{ m\in M: h(m)=t\}.
$$

\begin{definition} The function $h(m)$ is called an {\it
exhausting function} on $M$ if
\begin{enumerate}
\item
\label{h1}
$B_h (t)$ is precompact for all $t \in (0;h_0)$;
\item
\label{h2}
$\lim_{k \to \infty} h(m_k) = h_0$ along the every sequence $m_k \in M$
without accumulation points in $M\cup \partial M$;
\item
\label{h3}
$|\nabla h(m)| > 0 $ almost everywhere in $M$.
\end{enumerate}
\end{definition}

\begin{example} Let $M$ be a complete Riemannian manifold.
Then the distance function $h(m) = {\rm dist}(m_0,m)$ is an
exhausting function on $M$. Furthermore, at every regular point
of $h(m)$ we have $|\nabla h(m)| = 1$.
\end{example}

\smallskip
The following assertion provides a big family of exhausting functions
on minimal submanifolds in $\R{n}$.

\begin{lemma} %%%%%%%%%% Lemma 1.1.
\label{lem11} Let ${\mathcal M} = (M,u)$ be a minimal surface in
$\R{n}$ and $\varphi(x)$ be a smooth function defined on $\mathbb{R}^n$. We assume that the restriction $h(m) = \varphi\circ u(m)$
satisfies {\rm (1)} and {\rm (2)} in Definition {\rm 3} and there
exists a smooth positive function $\psi (t)$ for which the
composition $\psi\circ\varphi $ is a strong convex function. Then
$h(m)$ is an exhausting function on $M$.
\end{lemma}

\begin{proof}It suffices only to prove that $|\nabla h(m)|
> 0$ almost everywhere in $M$. If this fail then there exists a
closed set ${\mathcal Z}(h)$ of positive Lebesgues measure such
that $|\nabla h(m)| = 0$ if $m\in {\mathcal Z}(h)$. This would
implied existence of a point $m_0\in {\mathcal Z}(h)$ such that
the contingency of ${\mathcal Z}(h)$ at $m_0$ coincides with
$T_{m_0} M$. Let $\{E_k(m): \; k=1,\ldots,p\}$ be any orthonormal
system which is well defined in some neighbourhood of $m_0$. In
view of degeneracy of the gradient tangent component $\nabla h(m)
= (\grnd{}\varphi\circ u)^\top=0$ everywhere on ${\mathcal Z}(h)$,
we have for its normal component
$$
(\grnd{}\varphi\circ u)^\bot = \grnd{}(\varphi\circ u),\quad m \in
{\mathcal Z}(h).
$$
By the definition of the Weingarten map we notice that
\begin{equation*}
\begin{split}
A^\zeta (E_j) &= - (\grnd_{E_j}\grnd(\psi\circ\varphi))^\top = -
\left[\sum_{k=1}^{n}\grnd_{E_j}
\left(e_k\frac{\partial(\psi\circ\varphi)}{\partial x_k}
\right)\right]^\top\\
&=
-\sum_{k=1}^{n}\scal{E_j}{e_i^{\top}}e_k^{\top}
\frac{\partial^2(\psi\circ\varphi)}{\partial x_k \partial x_j},
\end{split}
\end{equation*}
where $\zeta = [\grnd{} (\psi\circ\varphi)]^{\bot} $ and
$\{x_k\}$ are the coordinate functions in $\R{n}$ which are
dual to the basis $\{e_k\}$. Taking the trace in the last
formula yields at $m_0$
$$
\scal{H}{\grnd(\psi\circ\varphi)} = -\sum_{i,k=1}^{n}\sum_{j=1}^{n}
\scal{E_j}{e_i^{\top}}
\scal{E_j}{e_k^{\top}}\frac{\partial^2(\psi\circ\varphi)}
{\partial x_k \partial x_j} =
$$
$$
 =-\sum_{i,k=1}^{n}\scal{e_k^{\top}}{e_i^{\top}}
\frac{\partial^2(\psi\circ\varphi)}{\partial x_k \partial x_i}.
$$
On the other hand, the Hessian
$$
\left|{\left|\frac{\partial^2(\psi\circ\varphi)}{\partial x_k
\partial x_i}\right|} \right|
$$
is a positively definite quadratic form at $m_0$ by virtue of
strong convexity of $\psi\circ\varphi $. Taking into
consideration vanishing of $ H \equiv 0$ and the last identity
we arrive at
$$
\scal{e_k^{\top}}{e_i^{\top}} = 0,
$$
for all $i,k \leq n$. Substituting $i = k$ gives $|e_k^\top|=0$, hence after
summation we obtain
$$
0 = \sum_{k=1}^{n}| e_k^{\top}|^2 = \dim M = p.
$$
The contradiction proves the lemma.
\end{proof}

\begin{example}
Let $(M;u)$ be a properly immersed surface in $\R{n}$. Then
$\varphi (m) = |u(m)|$ is an exhausting function on $M$.
\end{example}

In what follows we shall assume that any exhausting function $h$
satisfies the boundary condition (\ref{Neimann}) if the boundary
$\partial M$ is not empty.

\subsection{}Let $f(m)$ be a different from a constant function which
satisfies the maximum principle, i.e. for every open set $U \subset
M$ with compact closure
\begin{equation}
\max_{m \in \partial U}f(m) =
\max_{m \in \overline{U}}f(m).
\label{ue}\end{equation}

We notice that (\ref{ue}) valid for arbitrary $\p $-subharmonic
function. Really, it follows from the general maximum principle
for the elliptic type PDE's (see \cite{Trud}).

\begin{definition} A family of domains $\{D(\tau ):\tau
\in (\p,\beta )\}$ on $M$ is said to be the {\it asymptotic
tract} of $f(m)$ if:

\begin{enumerate}
\item[(i)]
for every $\tau \in (\p,\beta )$ the domain $D(\tau )$ is not empty
open component of the superlevel set $\{m \in M: f(m) > \tau\}$;

\item[(ii)]
$D(\tau_1) \supset D(\tau_2)$
for all $\tau_1 < \tau_2$ from the interval $(\p,\beta) $;

\item[(iii)]
either $\beta = +\infty $, or for some $\tau \in (\p,\beta )$
the following set is empty:
$$
D(\tau ) \cap \{m \in M : f(m) > \beta \}=\emptyset.
$$
\end{enumerate}
\end{definition}

It follows from (\ref{ue}), (i) and (iii) that any
$D(\tau )$ has a noncompact closure.

We say that two asymptotic tracts are {\it different} if
there exist  $\tau_1 \in (\p_1,\beta_1)$ and $\tau_2 \in
(\p_2,\beta_2) $ such that $D(\tau_1) \cap D(\tau_2) = \emptyset $.

Let $h(m)$ be an exhausting function on $M$ and $\{D(\tau ):\tau
\in (\p,\beta )\}$ be an asymptotic tract of $f(m)$.

\smallskip
{\bf Definition 5.} A tract $\{D(\tau )\}$ is called {\it
regular} if there exist $\tau_0 \in (\p,\beta )$ and a system
of pairwise disjoint intervals $\treu{k} \subset (0, h_0)$
converging to $h_0$ (i.e. $h_0\in\cup_{k} \overline{\treu{k}}$) such
that the set
$$
D(\tau_0) \cap \Sigma_h (t)
$$
contains no cycles for all  $t \in \cup_{k}
\Delta_k$. Otherwise, we say that $\{D(\tau )\}$ is a {\it
singular} asymptotic tract.

%%%%%%%%%%%%%%%%%%%%%%%%%%%%%%%%%%%%%%%%%%%%%%%%%%%%%%%%%%%%%%%%%%%%%%
\section{The fundamental frequency, its $N$-means and the Dirichlet
integral}
%%%%%%%%%%%%%%%\setcounter{equation}{0}
%%%%%%%%%%%%%%%%%%%%%%%%%% 2nd Section %%%%%%%%%%%%%%%%%%%%%%%%%%%%%%%

\subsection{} Let $\Sigma$ be a compact $p$-dimensional manifold. Further we
define the weighted fundamental frequency for the subsets of
$\Sigma$ which consist of a finite system of non-overlapping
$(p-1)$-dimensional compact submanifolds $\OO_j$ of $\Sigma$
with or without a boundary. Such a set we call to be {\it
simple} set (of $\Sigma$). We notice, that the $h$-sphere
$\Sigma_h(t)$ is a simple set for every regular value $t$ of
$h(m)$.

Let ${\U}$ be a simple set consisting of
components $\OO_{1},\OO_{2},\ldots,\OO_{k}$ and $\theta $ be a
positive in essential smooth function on ${\U}$. We say that a
Lipschitz function $\varphi $ {\it is admissible for} ${\U}$, or
$\varphi\wedge {\U}$, if for each component $\OO_{j}$ with
nonempty boundary $\partial\OO_{j}$:
$\varphi|_{\partial\OO_j}=0$ and the equality
$$
\int^{}_{\OO_{j}}|\varphi (m)|^{\p -2}\varphi (m)\theta (m) = 0,
$$
holds if $\partial \OO_j=\emptyset$.
Further we say {\it open} and {\it cyclic} components for the
components $\OO_j$ with or without a boundary respectively.

In Section 3 we use the following helpful property. Let $\OO_j$
is a cyclic component and $\varphi(m)$ be an arbitrary Lipschitz
function on $\OO_j$. Then the improvement function
$\varphi(m)-\xi$ will be already an admissible function for
$\OO_j$ if $\xi\equiv \phi_{\p}(\varphi,\OO_{j})$ is the unique
root of the following equation
$$
\int^{}_{\OO_{j}}\mm{\xi-\varphi (m)}^{\p -2}(\xi-\varphi (m))
\theta (m) = 0,
$$
where the integration considered over the $(p-1)$-dimensional
Hausdorff measure on $\OO_{j}$. One can show that the
last integral is increasing function on $\xi $, and hence
the quantity $\phi _{\p}(\varphi,\OO_{j})$ is well defined. The
explicit expression for $\phi_{\p}(\varphi,\OO_{j})$ is known only
for $\p = 2$:
\begin{equation}
\phi _{2}(\varphi,\OO_{j}) = \frac{\int_{\OO_j} \varphi
(m)\theta (m)} {\int_{\OO_j} \theta (m)}.
\label{fi}
\end{equation}

\smallskip

{\bf Definition 6.} The quantity
\begin{equation}
\lambda _{\p,\theta}({\U}) = \inf_{\varphi\wedge {\U}}
\left [
\frac{\int_{{\U}}\mm{\nabla\varphi}^{\p}\theta^{-1}}
{\int_{{\U}}\mm{\varphi}^{\p}\theta^{\p -1}}
\right ]^{\frac{1}{\p}}
\label{lambda}
\end{equation}
is called a {\it fundamental $\p $-frequency} of a simple
set ${\U}$ with respect to the weight $\theta $. We also
use the {\it reduced} fundamental frequency defined by
$$
\lambda^\star_{\p,\theta}({\U})=
\lambda_{\p,\theta}({\U})
$$
if ${\U}$ hasn't any cyclic component and
$\lambda^\star_{\p,\theta}({\U})=0$ otherwise.

\begin{example} If $\theta \equiv 1$ on ${\U}$, $\p = 2$ and
${\U}$ is a domain with the piecewise smooth nonempty boundary,
then $\lambda^2_{2,\theta}({\U})$ is the first nontrivial eigenvalue
of the classical variational Dirichlet problem. Other words,
$$
\Delta _{\Sigma}f + \lambda_{2,\theta}^2({\U}) f = 0,
$$
everywhere in ${\U}$ for some positive function $f(m)$ equal to
zero on $\partial {\U}$. Here $\Delta _{\Sigma}$ is the Laplace
operator on $\Sigma$.
\end{example}

\begin{example} Let now $\Sigma = \Sigma(\tau)$ be the submanifold in
$M$ given by the equation $h(m) = \tau$. Then it turns out that the
most natural choice of weight function $\theta (m)$ is $\mm{\nabla
h(m)}$. We write $\lalfa{h}{\OO}$ instead of
$\lalfa{\mm{\nabla h}}{\OO}$.
\end{example}

Now we formulate some basic properties of the fundamental
$\p$-frequency $\lambda^{(\star)}({\U})\equiv
\lambda^{(\star)}_{\p,\theta}({\U})$.

\begin{lemma} %%%%%%%%% 2.1.
\label{lem21}
If ${\U}_{1} \subset {\U}_{2}$, then $\lambda^\star ({\U}_1) \geq
\lambda^\star ({\U}_{2})$, i.e. $\lambda^\star$ is non-increasing
function of a set. If ${\U}$ consists of open components
$\OO_{1},\ldots,\OO_{k}$ only then
$$
\lambda ({\U}) = \min_{1\leq i\leq k} \lambda (\OO_i).
$$
\end{lemma}

\noindent
{\bf Proof.} If ${\U}_1$ contains a cyclic component then ${\U}_2$
does too. It follows that $\lambda^\star({\U}_1)=\lambda^\star
({\U}_{2})=0$ in this case. On the other hand, if ${\U}_2$ is without
cycles then such does ${\U}_1$ too. Hence monotonicity
of $\lambda^\star $ follows from the fact that $\varphi \wedge
{\U}_1$ implies $\varphi \wedge {\U}_2$.

To prove the second property we fix arbitrary functions $\varphi
_{i}\wedge\OO_{i}$ and set $\varphi (m) = \varphi _{i}(m)$
for $m \in \OO_i$. Hence
$$
\lambda ^{\p}(\OO_{i}) \int_{\OO_{i}} \mm{\varphi}^{\p} \leq
\int_{\OO_{i}} \mm{\nabla\varphi}^{\p}
$$
Summing of all this inequalities
$$
\min_{1\leq i \leq k} \lambda ^{\p}(\OO_{i})\,
\sum_{i=1}^{k}\int_{\OO_{i}} |\varphi|^{\p}\leq
\sum_{i=1}^{k}\int_{\OO_{i}} \mm{\nabla\varphi }^{\p}
$$
we obtain $\min_{1\leq i\leq k}\lambda (\OO_{i}) \leq \lambda (\OO
)$. The reverse inequality follows now from the first assertion
of the lemma.

The following assertion is an extension of the known for case $\p=2$
lower estimation \cite{C-Y} of the fundamental frequency of an open set.
Our proof is different from \cite{C-Y} where the main tool is the
maximum principle for solutions of the elliptic type PDE's.

\begin{lemma}%%%%%%%%%%%% 2.3.
\label{lem23}
Let $\p > 1$ and $f(m)>0$ be a function of $C^{2}({\U})$ such
that $\Z(f)=\emptyset$.
Then
\begin{equation}
\left ( \lalfa{\theta}{{\U}}\right )^{\p} \geq
\inf_{\U}\left [ - \frac{1}{(f \theta)^{\p -1}} {\rm
div} \frac{\mm{\nabla f}^{\p -2} \nabla f}{\theta} \right ].
\label{Yau}
\end{equation}
Here $\nabla$, ${\rm div}$ are considered with respect to the
inner metric of $\U$.
\end{lemma}

\noindent
{\bf Proof.} We fix a function $\varphi (m)\wedge {\U}$ and
denote by $\beta $ the right side of (\ref{Yau}).
Then from positivity of $f(m)$ we have everywhere in ${\U}$
$$
{\rm div }\left(\mm{\nabla f}^{\p -2} \frac{\nabla
f}{\theta}\right) + \beta (f\theta )^{\p -1} \leq 0,
$$
and hence,
\begin{equation}\label{esdiv}
\begin{split}
&\theta^{\p -1} \varphi^{\p}\beta =
\frac{\varphi^{\p} \beta}{f^{\p -1}}(f \theta)^{\p -1}
\leq - \frac{\varphi^{\p}}{f^{\p -1}} {\rm div}
\left(\mm{\nabla f}^{\p -2} \frac{\nabla f}{\theta}\right)\\
&=
\frac{\mm{\gd{f}}^{\p -2}}{\theta f^{\p}}\scal{\gd{f}}
{\p f\varphi^{\p -1}\gd{\varphi} - (\p -1)\varphi^{\p}\gd{f}}
-{\rm div}\left (\frac{\mm{\nabla f}^{\p -2}\varphi^{\p}}
{\theta f^{\p -1}} \nabla f\right )
\\
&\leq
\p\frac{(\varphi\mm{\nabla f})^{\p -1}}
{\theta f^{\p -1}} \mm{\nabla \varphi} - (\p
-1)\frac{\mm{\gd{f}}^{\p}\varphi^{\p}}{\theta f^{\p}}
- {\rm div}\left (\frac{\mm{\nabla f}^{\p -2}\varphi^{\p}}
{\theta f^{\p -1}} \nabla f\right ).
\end{split}
\end{equation}

Let ${\U}$ be decomposed into collection of the components
${\U}_{1},\ldots,{\U}_{k}$ with nonempty boundaries and the closed
components $\OO_{1},\ldots,\OO_{s}$. By virtue of the $\varphi
(m)\wedge {\U}$ we have $\varphi |_{\partial {\U}_i}=0$. Using
the Stokes` formula and (\ref{esdiv}) we arrive at the
inequality
$$
\beta \int_{{\U}} \theta ^{\p -1}\varphi^{\p} \leq \p
\int_{{\U}} \frac{(\varphi \mm{\gd{f}})^{\p -1}}{\theta f^{\p
-1}}\mm{\gd{\varphi}} \, - \, (\p -1)
\int_{{\U}}\frac{\mm{\gd{f}}^{\p}\varphi^{\p}}{\theta f^{\p}}.
$$
Applying the Cauchy`s integral inequality we obtain
$$
\beta \int_{{\U}} \theta ^{\p -1}\varphi ^{\p} \leq
\p \left [\int_{{\U}} \frac{\mm{\gd{f}}^{\p}\varphi^{\p}}
{\theta f^{\p}}\right ]^{\frac{\p -1}{\p}}
\left [\int_{{\U}}\mm{\gd{\varphi}}^{\p}{1\over \theta }\right ]^{1\over \p}-
$$
\begin{equation}
(\p -1) \int_{{\U}}\frac{\mm{\gd{f}}^{\p}\varphi^{\p}}{\theta f^{\p}}.
\label{esint}
\end{equation}
Given $k>1$, we can find $\varphi(m)\wedge {\U}$ such that
$$
\int_{{\U}}\frac{\mm{\gd{\varphi}}^{\p}}{\theta } \leq k \lambda^{\p}({\U})
\int_{{\U}} \varphi^{\p} \theta^{\p -1}.
$$
where we abbreviate $\lambda ({\U}) = \lalfa{\theta}{{\U}}$. It follows from
(\ref{esint})
$$
\beta \int_{{\U}} \varphi^{\p} \theta^{\p -1} \leq \p\lambda
({\U})k^{\frac{1}{\p}} \left
[\int_{{\U}}\frac{\mm{\gd{f}}^{\p}\varphi^{\p}}
{\theta f^{\p}}\right ]^{\frac{\p -1}{\p}}
\left [\int_{{\U}} \varphi^{\p} \theta^{\p -1}\right
]^{\frac{1}{\p}} -
$$

$$
(\p -1) \int_{{\U}}\frac{\mm{\gd{f}}^{\p}\varphi^{\p}}
{\theta f^{\p}}.
$$
Hence we obtain
\begin{equation}
A(\xi ) \equiv \beta \xi^{\p} - \p k^{{1\over \p}}
\lambda ({\U})\xi + (\p - 1) \leq 0,
\label{estA}
\end{equation}
where
$$
\xi =\left [ \int_{{\U}}\frac{\mm{\gd{f}}^{\p}\varphi^{\p}}{\theta
f^{\p}} \right ] ^{-\frac{1}{\p}}
\left [\int_{{\U}} \varphi^{\p} \theta^{\p -1}\right ]^{\frac{1}{\p}}.
$$

Now finding by standard arguments the minimum of $A(\xi )$ over all
$\xi \geq 0$ we conclude that
$$
A(\xi ) \geq (\p -1)\left [\frac{\beta}{\lambda ({\U})k^{\frac{1}{\p}}}
\right ]^{{1\over \p -1}} \; - (\p - 1)k^{{1\over \p}} \lambda ({\U})
$$
and by virtue of (\ref{estA}) we obtain
$$
\beta \;\leq \;\lambda^{\p}({\U}) k.
$$
Taking into account the condition on $k$ we complete the proof of
Lemma~\ref{lem23}.

\subsection{} We consider a $(p-1)$-dimensional subset $\OO \subset
\Sigma$ and an integer $N \geq 1$. Let us introduce
$$
\lambda (\OO, N) \equiv \lambda_{\alpha,\theta}(\OO
;N)=\inf\frac{1}{N} \sum^{N}_{i=1}\reduc_{\p,\theta}(\OO_{i}),
$$
where the infimum is taken over all systems
$\{\OO_{i}\}_{i=1}^{N}$ consisting of $N$ pairwise
non-overlapping simple subsets $\OO_{i} \subset \OO$.

\begin{definition} The quantity $\lambda (\OO, N)$ is called
$N$-{\it mean} of the fundamental frequency of $\OO$.
\end{definition}

\begin{remark}
It is clear from the definition that $\lambda(\OO;N)=0$ if $\OO$
contains $N$ or more cyclic components. Moreover, we give
without proof the following elementary property of this
characteristic (see \cite{Mik1}),
$$
\lalfa{\theta}{\OO ; 1} \leq
\ldots \leq \lalfa{\theta}{\OO ; N} \leq \ldots
$$
\end{remark}

\begin{lemma} %%%%%%%%% 2.2.
\label{lem22}
Let $\OO$ be an open component, $\theta \equiv 1$ and
$\varphi $ be an admissible for $\OO$ function which satisfies
equation $\Delta _{\p}\varphi=-\mu^{\p}\varphi \mm{\varphi }^{\p
- 2} $ in $\OO$. Then
$$
\lambda_{\p,\theta}(\OO, N) \leq \mu,
$$
where $N$ is the number components of the set $\OO_{0} =
 \{ m \in \OO : \mm{\varphi (m)} \ne 0\}$.
\end{lemma}

\noindent
{\bf Proof}. The inequality follows immediately from the
inequality $\lambda _{\p}(\OO_i) \leq \mu $ which is true for
each component $\OO_{j} \subset \OO_{0}$.

\subsection{} In this paragraph we give some estimates of the
fundamental frequency of one-dimensional sets (families of
curves) and their $N$-means. Let us consider a compact manifold
$\Sigma$, $\dim \Sigma = 1$. Then we have rather complete
information about $N$-mean of the fundamental frequency.

\begin{lemma} %%%%%% 2.4.
\label{lem24}
Given an open subset $\OO = \cup_{i=1}^{s}\OO_{i}
\subset \Sigma$ such that $\partial\OO_i \not\equiv \emptyset,
\OO_i\,\cap\,\OO_j = \emptyset$ and for any integer $N$ we have
\begin{equation}
\lambda_{2,\theta}(\OO) = \pi \left [ \max_{1\leq i\leq s}
\int_{\OO_i}\theta (t) \;dt
\right ]^{-1} \geq \frac{\pi}{\int_{\OO}\theta (t) \;dt}
\label{estlam}
\end{equation}
 and
\begin{equation}
\lambda_{2,\theta}(\OO;N) \geq\pi N\left[\max_{1\leq i\leq s}
\int_{\OO}\theta (t) dt \right ]^{-1}
\label{estlamN}
\end{equation}
\end{lemma}

\noindent
{\bf Proof.} Without loss of generality we can assume that $\OO$
is disjoint collection of intervals $\Delta _{1},\Delta
_{2},\ldots ,\Delta _{s}$ which are equipped by the natural
parametrization with their lengths and let $\mm{\Delta _{i}}$
be the length of $\Delta _{i}$. We choose arbitrary $\varphi
_{i}(t)\wedge \Delta _{i}$ and for $0 \leq t \leq \mm{\Delta
_{i}}$ set
$$
\xi _{i}(t) = \int\limits^{t}_{0} \theta _{i}(\tau)d\tau \,,\quad
$$
where $\theta _{i}(\tau)$ is the restriction of $\theta (t)$ on
$\Delta _{i}$. Applying the Wirtinger`s inequality we obtain

$$
\frac
{\int\limits_{0}^{\mm{\Delta_i}} \mm{\varphi_i '(t)}^2 \theta_i^{-1} dt}
{\int\limits_{0}^{\mm{\Delta_i}} \varphi_i (t)^2 \theta_i dt} =
\frac
{\int\limits_{0}^{\xi_i} \mm{\psi_i '(t)}^2 dt}
{\int\limits_{0}^{\xi_i} \psi_i (t)^2 dt} \geq
\left (\frac{\pi}{\xi _i}\right )^2.
$$

Here $\xi _{i} = \xi _{i}(\mm{\Delta _{i}})$ and
$\psi_i (\xi ) = \varphi_i\circ t(\xi )$. The equality holds
for $\psi_i (\xi ) = \sin (\xi\pi / \xi_i)$ and hence by virtue
of Lemma~\ref{lem21}, (\ref{estlam}) is proved.

In order to prove (\ref{estlamN}) we set $\epsilon > 0$ and
choose a decomposition $\OO$ into the collection of the
disjoint intervals $\delta _{i} \subset \OO$ such that
$$
\epsilon + \lambda_{2,\theta}(\OO ; N) \geq {1\over N}
\sum_{k=1}^{N} \lambda _{2,\theta}(\delta _{i}) =
{\pi \over N} \sum_{k=1}^{N}\left (\int_{\delta _{i}}\theta
_{i}(t)\; dt \right )^{-1}.
$$
Applying the inequality between the arithmetic and harmonic means
we obtain
$$
\epsilon + \lambda _{2,\theta}(\OO ; N) \geq
\pi N \left [{1\over N} \sum_{k=1}^{N} N\int_{\delta_{i}}\theta
_{i}(t) dt \right ]^{-1}
\geq \pi N \left [ \int_{\OO}\theta (t) dt \right ]^{-1},
$$
and in view of arbitrariness of $\epsilon > 0$ we arrive at the
required inequality.

The first part of the previous proof can be carried through
unchanged for the cycles. Hence, we have

\begin{lemma} %%%%% 2.5.
\label{lem25}
In the notations of the previous lemma and if
$\partial\OO_{i} = \emptyset$ for all $1 \leq i \leq s$,
\begin{equation}
\lambda _{2,\theta}(\OO) = 2\pi \left [ \max_{1\leq i\leq s}
\int_{\OO_{i}}\theta (t) dt
\right ]^{-1}.
\label{estLam}
\end{equation}
\end{lemma}

As for the higherdimensional case we have only fragmentary facts.
Let us assume that $M=\R{n}$ and $h(x)$ =$\mm{x}$ is the Euclidean
distance from the origin. It is clear that $\mm{\gd{h(x)}} = 1$ and
$\lalfa{\theta}{\OO}$ is equal to the ordinary fundamental frequency
$\lambda _{\p}(\OO)$. Let $n \geq 2$, then $h$-sphere is the Euclidean
sphere $S^{n-1}(t)$ of radius $t$, and given open subset $\OO
\subset S^{n-1}(t)$, we have (see, e.g. \cite{Mik1}),

\begin{equation}
\lambda (\OO,N) \geq
\left\{
\begin{array}{lcr}
c_{1}(n)\bigl (\frac{N}{{\rm meas}\; \OO}\bigr )^{\frac{1}{n-1}}&\mbox{if}
 &{\rm meas}\; \OO \leq \frac{1}{2}\omega_n t^{n-1},\\
&&\\
\frac{N-1}{N}c_{1}(n)\bigl (\frac{N}{{\rm meas}\; \OO}\bigr )^{\frac{1}{n-1}}
&\mbox{if} &{\rm meas}\; \OO \geq \frac{1}{2}\omega_n t^{n-1},
\end{array}
\right.
\label{miklyukov}
\end{equation}
where ${\rm meas}\; \OO$ is $(n - 1)$-dimensional Lebesgues measure
on $S^{n-1}(t)$, $\omega _{n}$ = ${\rm meas}\; S^{n-1}(1)$ and
$c_{1}(n)$ depends only on n.

\subsection{} Let $\D \subset M$ be an open subset and $f \in C^{1}(\D )$. For
the given $\p > 1$ we introduce the {\it Dirichlet integral} as
$$
\J{\p}{f}{{\mathcal D}} \equiv\int_{\D} \mm{\gd{f}}^{\p}.
$$
if the last integral is different from the infinity.

\begin{definition} Let $P,Q \subset \D$ be disjoint closed
sets. We define a {\it capacity} of the capacitor $(P,Q ;\D )$ by
\begin{equation}
{\rm cap}_{\p}(P,Q ;\D ) = \inf \int_{\D}\mm{\gd{\varphi}}^{\p}
\label{capac}
\end{equation}
where the infimum is taken over all locally Lipschitz functions
$\varphi : M \rightarrow \R{1}$ such that $\varphi (m) = 1$ when
$m \in P$ and $\varphi (m) = 0$ when $m \in Q$.
\end{definition}

We say that a manifold $M$ has an $\p $-{\it parabolic type} if
for any compact $F \subset M$ there is a sequence $\D_{1}\subset
\D_{2}\subset \ldots \subset \D_{k}\subset\ldots,\;
\cup_{k=1}^{\infty} \D_k = M$ of open sets $\D_{k}
\supset F$ with compact closures such that
$$
%%\begin{equation}
\lim_{k\to\infty}{\rm cap}_{\p }(F,M\setminus\D_{k}; M) = 0
%%\label{parab}
%%\end{equation}
$$
It is easy to see that the last limit doesn't`t depend on the
exhaustion sequence $\{\D_k\}$. Furthermore, we can assume that the family
$\{\D_{k}\}$ is a sequence of $h$-balls for the exhausting
function $h(m)$.
\par
\par
Given an $\p $-subharmonic function $f(m)$ on $M$, let us denote
by $\D = \D(\tau)$ a fixed element of the asymptotic tract $\{\D (\tau
)\}$ of $f(m)$. Let $t \in (h(\D );h_{0})$, where
$$
h(\D ) = \inf
\{h(m): m \in \D \},
$$
and
$$
P(t) = \overline{\D} \cap \overline{B_h(t)},
Q(t) = \overline{\D} \setminus {B_h(t)}.
$$

\begin{lemma} %%%%%%%% 2.6.
\label{lem26}
Let $M(t) = \max_{\D\cap \Sigma_h (t)} \{\tau,
f(m)\}$. Then for any numbers $t_{1}<t_{2}$ from the interval
$(h(\D ),h_{0})$:
\begin{equation}
\int_{P(t_{1})} \mm{\gd{f}}^{\p} \leq \p^{\p} \capac
(P(t_{1}),Q(t_{2});\D ) M^{\p}(t_{2}).
\label{estdirich}
\end{equation}

\end{lemma}

\noindent
{\bf Proof}.
Let $\varphi $ be any function which is admissible for the
calculation of $(P(t_{1}),Q(t_{2});\D )$. Then
$f_{1}(m)=(f(m)-\tau)$ vanishes everywhere in the set
$\partial\D \cap ({\rm int}M)$ and $\scal{f_{1}(m)}{ \nu} = 0$
on $\partial M$. Applying Stokes` formula we obtain
$$
\int_{\D} \varphi^{\p} \scal{\gd{f_1}}{\gd{f_1}\mm{\gd{f_1}}^{\p -2}}
\leq - \p \int_{\D} \varphi ^{\p -1} f_{1}
\scal{\gd{f_{1}}\mm{\gd{f_1}}^{\p -2}}{\gd{\varphi}}.
$$
Hence we have
\begin{equation}
\int_{\D} \varphi ^{\p}\mm{\gd{f}}^{\p} \leq
\p \left ( \int_{\D} \varphi ^{\p}\mm{\gd{f}}^\p \right )
^{\frac{\p -1}{\p}}\left (\int_{\D} \mm{f_1}^\p \mm{\gd{\varphi}}^{\p}
\right )^{\frac{1}{\p}}.
\label{ocen}
\end{equation}
and the maximum principle for $f_{1}(m)$ yields
$$
\left ( \int_{\D} \varphi ^{\p}\mm{\gd{f}}^\p \right )
\leq \p^{\p} \max_{\D\cap\Sigma_h(t_2)}f_1^{\p}(m)\;
\int_{\D\cap B_h(t_2)} \mm{\gd{\varphi}}^\p.
$$
Applying now the equality $\varphi (m) \equiv 1$ on $P(t_{1})$ and
taking infimum over all $\varphi (m)$ we obtain the
required estimate and Lemma~\ref{lem26} is proved.

\begin{lemma} %%%%%%%% 2.7.
\label{lem27}
Let $f(m)$ be an $\p $-harmonic function satisfying
(\ref{Neimann}) and $M(t) = \max
\{\mm{f(m)} : m \in \Sigma_h(t) \}$. Then for any
$t_{1}, t_{2}$ from $(0;h_{0})\;,\; t_{1} < t_{2}$ we have
$$
\int_{B_{h}(t_{1})} \mm{\gd{f}}^{\p} \leq \p^\p
\left [
\int\limits_{t_1}^{t_2}\frac{M^{-\frac{\p}{\p
-1}}(t)}{S(\Sigma_t)}\, dt
\right ]^{1-\p},
$$
where $S(\Sigma_t) = \left (\int_{\Sigma_{h}(t)} \mm{\gd{h}}^{\p
-1}\right ) ^{\frac{1}{\p -1}}$ is the flow of the
exhausting function.
\end{lemma}

The {\bf proof} differs only in the details from the same one given
above. Namely, from the coarea formula (\cite{Fed}, Theorem 3.2.22)
and (\ref{ocen}) it follows
\begin{equation}
\int_{B_{h}(t_{1})} \mm{\gd{f}}^{\p} \leq \p^{\p}
\int_{B_h(t_2)} \mm{f}^\p \mm{\gd{\varphi}}^\p \leq
\p^{\p} \int\limits_{t_1}^{t_2} M^{\p}(t) dt \int_{\Sigma_h(t)}
\frac{\mm{\gd{\varphi}}^{\p}}{\mm{\gd{h}}}
\label{efm}
\end{equation}
for any admissible for the capacitor $(\overline{B_h(t_{1})},
M\setminus B_h(t_{2}); M)$ function
$\varphi (m)$
\par
Set
$$
G(t) = \int^{t}_{t_{1}} {d\xi \over M^{\p\over \p -1}(\xi
)S(\Sigma(\xi))}
$$
and substitute ${G\circ h(m)/G(t_{2})}$ instead of $\varphi (m)$ into
(\ref{efm}). Hence, we obtain
$$
\int_{B_{h}(t_{1})} \mm{\gd{f}}^\p \leq \p^{\p} G(t_{2})^{1-\p},
$$
as required.

\subsection{} Let us suppose that $h(m)$ be a positive $\p $-harmonic
function everywhere outside of its own null-level set (the main
case is $h$ to be the absolute value of a $\p$-harmonic function).
Then the flow $S(\Sigma_t)$ defined in Lemma~\ref{lem27} does
not depend on $t$. The quantity
$$
S(h)\equiv S(\Sigma_{t}) = \int_{\Sigma_{t}(h)}
\scal{\nu }{\gd{h}\mm{\gd{h}}^{\p -2}},
$$
is said to be a {\it full flow} of $h$.

\begin{lemma} %%%%%% 2.8.
\label{lem28}
 For any $t_{1},t_{2}$ from $(0,h_{0})$ such that $t_{1} <
t_{2}$,
$$
\capac(\overline{B_{h}(t_{1})}, M\setminus B_{h}(t_{2}); M)=
\left ({S(h)\over t_{2} - t_{1}}\right )^{\p -1}.
$$

\end{lemma}

\noindent
{\bf Proof}. We notice that the function $\varphi (m)$ defined
for $m\in \overline{B_h(t_2)} \setminus B_{h}(t_{1})$ as
$$
\varphi (m) = {t_{2} - h(m)\over t_{2} - t_{1}},
$$
equal to $0$ when $m \in M \setminus B_{h}(t_{2})$ and $1$ when
$m \in B_{h}(t_{1})$ is admissible for the capacitor
$(\overline{B_{h}(t_{1})}, M\setminus B_{h}(t_{2}); M)$. Then
from (\ref{capac})
we have
$$
{\rm cap}_{\p,h}(t_{1},t_{2}) \leq
\int_{B_h(t_2)\setminus\overline{B_h(t_1)}}
\mm{\gd{\varphi}}^{\p},
$$
where we set
$$
{\rm cap}_{\p,h}(t_{1},t_{2}) =
\capac(\overline{B_{h}(t_{1})}, M\setminus B_{h}(t_{2}); M).
$$
Applying now $\p $-harmonicity of $h(m)$ and the Stokes'
formula we obtain
$$
\int_{B_h(t_2)\setminus\overline{B_h(t_1)}}\mm{\gd{\varphi}}^{\p} =
\int_{B_h(t_2)\setminus\overline{B_h(t_1)}} ( \mm{\gd{\varphi}}^{\p} +
\varphi \;{\rm div} \mm{\nabla\varphi}^{\p-2}\gd{\varphi} ) =
$$

$$
=\int_{\Sigma_h(t_2)} \varphi\mm{\gd{\varphi}}^{\p-1}\scal{\gd{\varphi
}}{\nu} -
\int_{\Sigma_h(t_1)} \varphi\mm{\gd{\varphi}}^{\p-1}\scal{\gd{\varphi
}}{\nu}.
$$
Now taking into account the boundary conditions on
$\varphi (m)$ we arrive at
$$
{\rm cap}_{\p,h}(t_{1},t_{2}) \leq
\left ({S(h)\over t_{2} - t_{1}}\right )^{\p -1}.
$$
To verify that in fact the equality holds, it is sufficient to
observe that $\varphi (m)$ is actually an extreme for the
calculation of the $\p $-capacity of the capacitor
$(P(t_{1}),Q(t_{2}))$, since it is $\p $-harmonic and satisfies
the "natural" boundary condition $\scal{\gd{h}}{\nu } = 0$ on
$\partial M$ (see e.g. \cite{C}).

%%%%%%%%%%%%%%%%%%%%%%%%%%%%%%%%%%%%%%%%%%%%%%%%%%%%%%%%%%%%%%%%%%%%%%%
\section{The Main Inequality}
\label{sec3}
%%%%%%%%%%%%%%%%%%%%%%%%% 3rd Section %%%%%%%%%%%%%%%%%%%%%%%%%%%%%%%%%

\subsection{} Let $M$ be a noncompact $p$-dimensional manifold with a
fixed exhausting function $h(m) : M \rightarrow [0,h_{0})$. We
assume that $f$ is an $\p $-subharmonic function satisfying
(\ref{Neimann}) and $\{\D(\tau)\}$ is an asymptotic tract of
{\it f}. Let us choose an arbitrary subset $\D\equiv\D
(\tau_{0})$, where $\tau_0$ is a regular value of both $h$ and $f$,
and set $f_{1}(m)=f(m)-\tau_{0}$ everywhere in $\D$. Let
$h(\D)=\inf_{m\in\D}h(m)$. Then for any regular value $t >h(\D
)$ of $h(m)$ we have
$$
J(t) \equiv \J{\p}{f}{\D\cap B_h(t)}= \int\limits_{\D\cap B_h(t)}|\nabla
f| ^{\p } \le \int\limits_{\D\cap \Sigma_{h}(t)}
f_{1}\scal{\nabla f_{1}| \nabla f_{1}|^{\p -2}}{{\nabla h\over |
\nabla h| }},
$$
because $f_{1}(m)\Delta _{\p }f_{1}(m) \ge 0$ everywhere in $\D$.
Here ${\nabla h/| \nabla h| }$ is the outward normal $\nu (m)$
to $\Sigma_{h}(t)\cap\D$ with respect to $M$. The set
$\Sigma_{h}(t) \cap \D $ is compact and, hence, it splits into
the cycles $\Gamma_{1}(t),\ldots,\Gamma _{l}(t)$ and components
$\gamma _{1}(t),\ldots ,\gamma _{k}(t)$ with nonempty
boundaries. Let $\Gamma (t) =\cup_{i=1}^{l} \Gamma_{i}(t)$ and
$\gamma (t)=\cup_{j=1}^{k}\gamma _{j}(t)$. We set
$$
\lambda (t)=\min \{\lambda_{\p,\nabla h}(\gamma_i),
\lambda_{\p,h}(\Gamma_j)\}
$$
and $q_{j}(t) = \phi _{\p }(f_{1}(m);\Gamma _{j}(t))$. Then
the function
\par
$$
\upsilon (m) = \left\{
\begin{array}{llr}
f_{1}(m),& m \in \gamma (t) & \\
f_{1}(m) - q_{j}(t),& m \in \Gamma _{j}(t)&\mbox{for}\;\; 1\leq
j\leq l
\end{array}
\right.
$$
is admissible for $\Sigma_h(t)\cap \D$. Applying the Minkowski
inequality we have
\par
$$
\upsilon \scal{\nabla f_{1}| \nabla f_{1}| ^{\p -2}}{\nabla h} \le {\p -
1\over \lambda (t)\p }
\left |
\scal{\nabla f_{1}| \nabla f_{1}
| ^{\p -2}}{{\nabla h\over | \nabla h| }}
\right |
^{{\p \over \p -1}}+
$$

$$
{1\over \p } \lambda ^{\p -1}(t)|\upsilon| ^{\p }| \nabla h| ^{\p },
$$
and, therefore,
\par
$$
J(t) \le \int\limits_{\Sigma_{h}(t)\cap \D }\upsilon \scal{\nabla f_{1}|
\nabla f_{1}
| ^{\p -2}}{{\nabla h\over | \nabla h| }} +
\sum^{l}_{j=1}q_{j}(t)\int\limits_{\Gamma _j} \scal{\nabla f_{1}|
\nabla f_{1}| ^{\p -2}}
{{\nabla h\over | \nabla h| }}\le
$$
$$
{\lambda ^{\p -1}(t)\over \p }
\int\limits^{}_{\Sigma_{h}(t)\cap\D }
|\upsilon| ^{\p }| \nabla h| ^{\p-1 } +
$$
$$
{\p - 1\over \p \lambda
(t)}\int\limits^{}_{\Sigma_{h}(t)\cap\D
}{1\over | \nabla h| }\left | \scal{\nabla f_{1}| \nabla f_{1}|
^{\p -2}}{{\nabla h\over | \nabla h| }} \right |^{{\p \over \p
-1}} + Q(t),
$$
where
\begin{equation}
Q(t) =\sum^{l}_{j=1}q_{j}(t) \int\limits^{}_{\Gamma
_{j}(t)}\scal{\nabla f_{1}| \nabla f_{1}| ^{\p -2}}{\nu}.
\label{Qut}
\end{equation}
Taking into account the definition of $\lambda(t)$ we arrive at
$$
J(t)\leq \frac{1}{\p\lambda(t)} \int\limits_{\Sigma_{h}(t)\cap\D }
\frac{|D\upsilon|^{\p }}{|\nabla h|}+
$$
\begin{equation}
\frac{\p-1}{\p\lambda (t)}
\int\limits_{\Sigma_{h}(t)\cap\D}\frac{|\nabla f_1|^{\p}}{|\nabla h|}
\left |\bigl\langle{\frac{\nabla f_{1}}{|\nabla f_1|}},{\frac{\nabla h}{|
\nabla h|}}\bigr\rangle \right |^{{\p\over \p-1}} +Q(t),
\label{Jit}
\end{equation}
where $D$ is the induced covariant derivative on $\Sigma (t)$.

By virtue of orthogonality of $\nu=\nabla h(m)|\nabla
h(m)|^{-1}$ to the tangent space $T_m\Sigma_{h}(t)$, we obtain
the following gradient decomposition formula
\par
\begin{equation}
\nabla f_{1} = D f_1 + \nu \scal{\nabla f_{1}}{\nu }.
\label{pifagor}
\end{equation}

Moreover, $D \upsilon = Df_{1}$ everywhere in $\D \cap
\Sigma_{h}(t)$.
Hence, from (\ref{Jit}), (\ref{Qut}) and
(\ref{pifagor}) we obtain
$$
J(t) \le Q(t) + {1\over \lambda (t)}
\int\limits^{}_{\Sigma_{h}(t)\cap\D}{| \nabla f_{1}|^\p \over |
\nabla h| }
\left ( {\xi ^{\p }\over \p } + {\p - 1\over \p
}(1-\xi^2)^{{\p\over 2(\p -1)}}\right ),
$$
where $\xi = | Df_{1}| /|\nabla f_{1}|\le 1$.
Analysis of the last expression in the parentheses shows that it
is decreasing function on $\xi$ for $\p\geq 2$ and increasing
one if $1<\p<2$. Hence,
$$
J(t)\leq Q(t) +
\frac{c(\p)}{\lambda(t)}\int\limits^{}_{\Sigma_{h}(t)\cap\D}{|
\nabla f_{1}|^\p
\over |\nabla h|},
$$
where
$$
c(\p)= \begin{cases} \displaystyle{\frac{\p-1}{\p}} &\mbox{for}
\p\geq 2;\cr \displaystyle{\frac{1}{\p}}&\mbox{for} \p\in
(1;2).\cr
\end{cases}
$$

Because the set $\Sigma_{h}(t)\cap \D$ is the $t$-level of $h$,
the co-area formula yields
\par
\begin{equation}
J(t) \le Q(t) + \frac{c(\p)}{\lambda (t)} {d J\over dt}
\label{JQA}
\end{equation}
for a.e. $t \in (h(\D ),h_{0})$.

\subsection{} In this paragraph we assume that there exists at least one
regular asymptotic tract of $f$. Let $\Delta $ be the open
subset of $(h(\D );h_{0})$ corresponding to regularity of
$\D(\tau_{0})$, i.e. $\Gamma(t)=\emptyset$ for $t\in\Delta$,
where
$$
\Delta = \bigcup_{i=1}^{\infty} \Delta _{i} \;,
\quad \Delta _{i} = (\p_{i}, \beta _{i}),
$$
and $h_0\in\overline{\Delta}$. We observe that the
characteristic $J(t)$ is nondecreasing function. Moreover, for
$t\in\Delta$ we have by Lemma \ref{lem21} that
$\lambda(t)=\reduc_{\p,h}(\Sigma_h(t)\cap \D)$, and
by virtue of (\ref{JQA}) we obtain
$$
J'(t) \geq c(\p)\reduc_{\p,h}(\Sigma_h(t)\cap \D)\, J(t)
$$
almost everywhere in $(h(\D),h_{0})$. By virtue of absolute
continuity of $J(t)$, we arrive at

\begin{theorem} %%%% 3.1
\label{th31}
Let $f(m)$ be an $\p$-subharmonic function on $M$ and $\D = \D (\tau_{0})$
be a regular component of an asymptotic tract $\{\D
(\tau)\}$. Then for any $t_{1}$, $t_{2}$ from $(h(\D
),h_{0})$
$$
\J{\p}{f}{\D \cap B_{h}(t_1)}\le
\J{\p}{f}{\D \cap B_{h}(t_2)} \exp \left (-
c(\p)\int\limits_{(t_{1},t_{2})\cap\Delta }
\lambda^\star_{\p,h}(\Sigma_{h}(t)\cap \D)dt \right ).
$$
\end{theorem}

\par
\begin{corollary} %%%%%% 3.1
\label{cor31}
Let $f$ be an $\p$-subharmonic function having at least $N$
different regular asymptotic tracts $\{\D_{i}(\tau)\},\ldots
,\{\D_{N}(\tau)\}$ and 
$$\sigma =\max_{1\leq i\leq N}
h(\D_{i}),$$ 
where $\D_{i} \equiv \D_{i}(\tau_{i})$ are regular
elements. Then for any $t_{1}< t_{2}$ from $(\sigma,h_{0})$
\begin{equation}
N \!\min_{1\leq i\leq N}\!\! \int\limits_{\D_i \cap B_{h}(t_1)}\!\! |\nabla f|^\alpha \le
\exp \biggl(- c(\p) \int\limits^{t_{2}}_{t_{1}} \reduc
_{\p,h}(\Sigma_{h}(t);N)\;dt
 \biggr)
\int\limits_{B_{h}(t_2)}\!\!|\nabla f|^\alpha
\label{Nmin}
\end{equation}
\end{corollary}

\noindent
{\bf Proof}. According to the previous theorem we have for every
domain $\D_{i}$ and for $t_{1} < t_{2}$, $t_i\in(\sigma,h_{0})$ that
$$
\exp \left ( c(\p) \int\limits^{t_{2}}_{t_{1}}
\lambda _{i}(t)dt \right )\;
\int\limits_{B_{h}(t_{1})\cap \D_i}|\nabla f| ^{\p }
\le \int\limits_{\D_{i}\cap B_h(t_2)}\mm{\gd{f}} ^{\p },
$$
where $\lambda _{i}(t) = \reduc _{\p,h}(\Sigma_{h}(t) \cap
\D_{i})$. Summing the last inequalities yields
\begin{equation}
\min_{1\leq i\leq N} \J{\p}{f}{\D_{i} \cap B_{h}(t_1)}\cdot
\sum^{N}_{i=1}\exp \left(c(\p)
\int\limits^{t_{2}}_{t_{1}} \lambda _{i}(t)\;dt \right)
\le
\J{\p}{f}{B_{h}(t_2)}%\mm{\nabla f} ^{\p }.
\label{sred}
\end{equation}

Applying the inequality between arithmetic and geometric means
in the left side of (\ref{sred}) we arrive at
$$
\min_{1\leq i\leq N} \J{\p}{f}{\D_{i} \cap B_{h}(t_1)}\cdot
N\exp\frac{1}{N}\sum^{N}_{i=1} \left(c(\p)
\int\limits^{t_{2}}_{t_{1}} \lambda _{i}(t)\;dt \right)
\le
\J{\p}{f}{B_{h}(t_2)}.
$$

But the domains $\D_{1},\ldots,\D_{N}$ are pairwise
non-overlapping ones, and therefore for every $t \in
(t_{1},t_{2})$
$$
{1\over N} \sum^{N}_{i=1} \lambda _{i}(t) = {1\over N}
\sum^{N}_{i=1}\reduc _{\p,h}(\D_i \cap \Sigma_{h}(t)) \ge
\lambda_{\p,h}(\Sigma_{h}(t),N).
$$
Hence, the required inequality is proved.

\begin{corollary}
\label{cor32}
Let $M$ be a manifold with an exhausting function $h(m)$ such
that for some integer $N \ge 1$ and $h_{1}<h_{0}$
\begin{equation}
\int\limits^{h_{0}}_{h_{1}} \lambda _{\p,h}(\Sigma_{h}(t);N)\,dt
= + \infty.
\label{ash}
\end{equation}
Then every $\p $-subharmonic function $f(m)$ with finite Dirichlet integral
$$
\int\limits^{}_{M} | \nabla f| ^{\p } < \infty
$$
have at most $(N - 1)$ different regular asymptotic tracts.
\end{corollary}

\noindent
{\bf Proof}. Indeed, if $f(m)$ has $N$ different regular
asymptotic tracts, then from (\ref{Nmin}) and (\ref{ash}) by
finiteness of the Dirichlet integral we conclude for some $i \le
N$ and $t_{i}\in(h(\D_{i}),h_{0})$ that everywhere in the open set
$\D_{i} \cap B_{h}(t_{i})$ one holds $| \nabla f| \equiv 0$. It
follows that $f(m) \equiv {\rm const}$ when $m \in \D_{i} \cap
B_{h}(t_{i})$, that contradicts the definition of
$\D_{i} = \D_{i}(\tau_{i})$.

\begin{corollary}
\label{cor33}
{\it Let $M$ be a manifold without a boundary and an
exhausting function $h$ satisfies
\par
\begin{equation}
\int\limits^{h_{0}}_{h_{1}} \lambda _{\p,h}(\Sigma_{h}(t);2)dt =
+\infty \label{hyper}
\end{equation}
Then constants are only $\p$-harmonic functions on $M$ with
finite Dirichlet integral}.
\end{corollary}

\noindent
{\bf Proof.} Let $\Delta=\cup_{i}\Delta_i$ be the subset of
$(h_1,h_0)$ consisting of the regular values $t$ of $h$ such
that $\lambda_{\p,h}(\Sigma_h(t),2)>0$. It follows from Remark 1
that $\Sigma_h(t)$ can contain one cyclic components at most. On
the other hand, $\Sigma_h(t_0)$ is a compact submanifold without
a boundary for all $t_0\in \Delta$. Hence, if $\Sigma_h(t_0)$
contains a cyclic component then it must contain two or more
such components. But this contradicts with the fact that $t_0\in\Delta$ and
shows hereby that $\Sigma_h(t_0)$ contains no cycles.

Now, let $f(m)$ be an arbitrary $\p$-harmonic function different
from a constant. Given a fixed point $m_{0}$ we consider the new
subharmonic function $f_{1}(m) = | f(m) - f(m_{0})|$ which has at
least two different asymptotic tracts. The last property
is a direct consequence of maximum principle for $f_1$. Then it
follows from the arguments above that these tracts are regular.

Other conditions concerning the Liouville type theorems
can be found, for instance, in \cite{Trud}, \cite{Gr}. We
emphasize that our assertions have been formulated without any
requirements on the geodesic completeness. The following example
shows that the conditions (\ref{ash}), (\ref{hyper}) do not
yield it.

\begin{example} Let $M$ be realized as the compact rotational
hypersurface in $\R{n}$ with a profile function $\rho
(x_{n})$, i.e.
\par
$$
x^{2}_{1} + x^{2}_{2} +\ldots + x^{2}_{n-1} = \rho
^{2}(x_{n}),\qquad x_{n} > 0.
$$
We suppose that $\rho (0) = 0$ and the original in $\R{n}$ is
the limit singular point of $M$. Namely,
$$
\rho (t) = \begin{cases} t &for t\in (0;1);\cr \zeta(t)&\mbox{for}
t\in (1;2);\cr \sqrt{9-t^2}& \mbox{for} t\in (2;3);\cr
\end{cases}
$$
where $\zeta(t)$ is a sufficient smooth sewing function.
Thus $M$ is a "drop" with regular top point
$A=(0,\ldots,0,3)$. Then the geodesic distance $h(m)$ from $A$ to
$m \in M$ is an exhausting function on $M$ and
\par
$$
h_{0} = \lim_{m\to 0} h(m) < + \infty.
$$
Therefore, $M$ isn`t geodesically complete. On the other hand, all
$h$-spheres $\Sigma_{h}(t)$ are $(n-2)$-dimensional Euclidean
spheres of radius $R(t)$. Here $R(t)=\rho (m)$ where $h(m) = t$.
It follows for the values of $h(m)$ sufficiently closed to $h_0$
that
$$
R(h(m)) =\frac{1}{\sqrt{2}} (h_{0} - h(m)).
$$
\end{example}

We observe also that $| \nabla h| \equiv 1$ everywhere in $M\backslash
\{A\}$. Hence, from (\ref{miklyukov}) we obtain
$$
\lambda _{\p,h}(\Sigma_{h}(t),N) = {c\over h_{0} - t},
\; \mbox{ as } \; t \to h_{0},
$$
where $c$ depends only on $n$, $N$. We see that (\ref{hyper})
is valid, while $M$ isn`t geodesic complete.
\par
\medskip

\subsection{} Now we assume that an $\p $-subharmonic function
$f(m)$ has at least $N \ge 1$ regular asymptotic tracts
$\{\D_{1}(\tau)\},\ldots ,\{\D_{N}(\tau)\}$ and choose
pairwise disjoint domains $\D_{k} = \D_{k}(\tau_{k}),
1 \le k \le N$ such that the sets $\D_{k} \cap \Sigma_{h}(t)$
would answer to the regularity's condition.
\par
Applying (\ref{estdirich}) and (\ref{sred}) we obtain
\begin{equation*}
\begin{split}
\min_{1\leq i\leq N} & \J{\p}{f}{\D_{i}\cap B_{h}(t_1)}
\sum^{N}_{k=1}\exp \left\{c(\p)\int\limits^{t_{2}}_{t_{1}}
\reduc_{\p,h}(\D_{k}\cap\Sigma_{h}(t))dt \right\}\\
&\le \p ^{\p } M^{\p }(t_{3}){\rm cap_{\p,h}} (t_{2}, t_{3}).
\end{split}
\end{equation*}
where
$$
M(t) = \max_{m\in \Sigma_h(t)} \{\tau_{1},\tau_{2},\ldots
,\tau_{N},f(m)\}$$
and the numbers $t_{1}< t_{2}< t_{3}$ belong
to $(h_{1},h_{0})$.

Arguing as in the proof of the corollary \ref{cor31} we conclude
$$
\min_{1\leq i\leq N}\J{\p}{f}{\D_{i}\cap B_{h}(t_1)}
\leq
$$
$$
\frac{\p ^{\p } M^{\p
}(t_{3})}{N} {\rm cap_{\p,h}} (t_{2}, t_{3})
\exp \left(-c(\p)\int\limits^{t_2}_{t_{1}}\lambda _{\p,h}
(\Sigma_{h}(t),N) \;dt\right).
$$

Thus, we have proved the following version of the Denjoy-Ahlfors
theorem.

\begin{theorem} %%%%%% 3.2
\label{th32}
Let $f(m)$ be an $\p $-subharmonic function and for some integer $N
\ge 1$
$$
%%\begin{equation}
\liminf_{t,\xi\to h_0} M(\xi)\,
{\rm cap}^{{1\over \p }}_{\p,h}(t,\xi )\,
\exp\left(-\frac{c(\p)}{\p}
\int\limits^{t}_{t_{1}}\lambda _{\p,h}
(\Sigma_{h}(t),N) dt\right)=0,
%%\label{denal}
%%\end{equation}
$$
where
$$
M(t)=\max_{\Sigma_h(t)} f^+(m),
$$
$t_{1} > 0$ is fixed and $t$, $\xi$ are approaching to
$h_0$ such that $t_1 < t < \xi$.
Then $f(m)$ has at most $(N-1)$ different regular asymptotic
tracts.
\end{theorem}

\section{The cyclic case: singular asymptotic tracts}
\label{sec4}

\subsection{} We observe (see for a details \cite{HK}, \cite{Mik1}) that the
topologically trivial manifolds such as $\R{n}$ or graphs over
$\R{n}$ admit a big store of subharmonic functions with regular
asymptotic tracts. On the other hand, there exist manifolds
having nontrivial harmonic functions with compact level-sets.
This means that for such a manifold the term
$\lambda(\Sigma_h(t);N)$ is trivial. It is a consequence of the
definition of this value by the reduced fundamental frequency.
To deal with such manifolds we must refine our tool. On the other
hand, there are some difficulties in the higher-dimensional case.
Namely, we need the lower estimations of the first eigenvalue of
compact manifolds. In this section we consider the two-dimensional
case only and will assume that $\p=2$.

\begin{definition} We say that a manifold $M$ admits the
{\it harmonic exhausting} if there exists an exhausting function
$h(m)$ of $M$ satisfying the property: $\Delta h(m) = 0$ everywhere
outside of $\Sigma_{h}(0)$. If $\partial M$ is not empty we assume
that the condition {\rm (\ref{Neimann})} holds for $h$ also.
\end{definition}

Let $f(m)$ be a harmonic function on $M$ and $t_{0}$ be a regular
value of $h(m)$. Then the characteristic $Q(t)$ introduced in
(\ref{Qut}) is differentiable at $t_{0}$. Really, it is a consequence
of the Stokes` formula and the property of orthogonality of $\nabla
h$ to $\Gamma _{j}(t)$ that
\begin{equation}\label{potok}
\begin{split}
{dq_j(t_{0})\over dt}& = {d\over dt}\biggl [
S^{-1}(h,\Gamma _{j}(t_{0}))\int\limits^{}_{\Gamma _{j}(t)}|
\nabla h| f
\biggr ]_{t=t_{0}} \\
&={1\over S(h,\Gamma _{j}(t_{0}))}\biggl [ {d\over dt}
\int\limits_{\Gamma_j(t)} f\scal{\nabla h}{\nu } \biggr
]_{t=t_{0}} = {1\over S(h,\Gamma _{j}(t_{0}))}
\int\limits_{\Gamma_j(t_0)}{\hbox{div}(f\nabla h)\over {| \nabla
h| }}\\
&={1\over S(h,\Gamma _{j}(t_0))} \int\limits_{\Gamma_j(t_0)}
\scal{\nabla f}{{\nabla h\over | \nabla h| } } = {S(f,\Gamma
_{j}(t_0))\over S(h,\Gamma _j(t_0))}.
\end{split}
\end{equation}
Here by $S(\varphi,\Gamma )$ we denote the flow
of $\varphi $ through $\Gamma $. By virtue of  (\ref{potok}) and
(\ref{Qut}) we obtain
$$
%%\begin{equation}
{dQ\over dt}(t_{0}) = \sum^{l}_{j=1} {S^{2}(f,\Gamma
_{j}(t_{0}))\over S(h,\Gamma _{j}(t_{0}))}.
%%\label{massa}
%%\end{equation}
$$
From harmonicity of $f$ and $h$ we obtain the following

\begin{proposition}
\label{prop41}
The derivative $\omega (t) = Q'(t)$ is a piecewise constant
function. Moreover, $Q(t)$ is a piecewise linear function.
\end{proposition}

We observe that the full flows of $f$ and $h$ through $\Gamma(t)$
does not depend on the value $t$. Moreover,
$$
S(h) = S(h,\Gamma (t))= \sum^{l}_{j=1}S(h,\Gamma _{j}(t))
=\int\limits^{}_{\Gamma (t)}| \nabla h| > 0
$$
and from Cauchy`s inequality we have
$$
%%\begin{equation}
\omega (t) = {dQ\over dt} \ge {S^{2}(f)\over S(h)}
%%\label{omg}
$$
%%\end{equation}

\subsection{} In this section we suppose that $h(m)$ has {\it isolated}
critical points only. We notice that this requirement is
realized when $\dim M = 2$. Really, due M.Morse, every harmonic
function $f(m)$ on $M$ is simultaneously a pseudoharmonic function in
local coordinates of $M$, and the needed assertion follows from the
standard properties of pseudoharmonic functions.

It follows that $Q(t)$ is continuous and hence
$$
Q(t) = \int\limits^{t}_{0} \omega (\eta )d\eta,
$$
Hence, the Main Inequality provides the following relation
\begin{equation}
{dJ(t)\over dt} \ge 2\lambda _{h}(\Sigma_{h}(t))\left ( J(t) -
\int\limits^{t}_{0} \omega (\eta )d\eta \right )
\label{pro}
\end{equation}
where
$$
J(t) = \int\limits_{B_h(t)} | \nabla f| ^{2}.
$$

We need the following consequence of Lemma~\ref{lem25}:

\begin{proposition}
\label{prop42}
 $M$ be a two-dimensional manifold and $h(m)$
be a harmonic exhausting on $M$. Then
\begin{equation}
\lambda _{h}(\Sigma_{h}(t)) = 2\pi \min_{1\leq j\leq l}
 \left\{S^{-1}(h,\Gamma _{j}(t))\right\} \ge 4\pi S^{-1}(h).
\label{Slam}
\end{equation}
\end{proposition}

\begin{proof} We notice by harmonicity of $h(m)$ that the set
$\Sigma_{h}(t)$ splits into the union of two nonempty level-sets
$\Sigma^{\pm}_{h}(t)=\{m: h(m) =\pm t \}$. Moreover, by the Stokes` formula we have
\par
$$
S(h,\Sigma^{+}_{h}(t)) = S(h,\Sigma^{-}_{h}(t)) =
{1\over 2} S(h),
$$
and positivity of the flow yields
\par
\begin{equation}
S(h,\Gamma _{j}(t)) \le S(h,\Sigma^{\pm}_{h}(t))
 = {1\over 2} S(h).
\label{prom}
\end{equation}
Hence (\ref{Slam}) follows from (\ref{estLam}) and (\ref{prom}).
\end{proof}

Let $e$ be a fixed unit vector in $\R{n}$ and $\{\Pi _{t}\}$ be
the family of all hyperplanes given by the equation $\scal{x}{e}
= t$. We notice that every $\Pi_t$ is orthogonal to $e$.

\begin{definition}[\cite{MV}]
A surface $\M= (M,x)$ is called {\it tubular} if there
exists an interval $(a;b)$ such that any portion of $\M$
situated between two hyperplanes $\Pi_{t_{1}}$ and $\Pi_{t_{2}}$
with $t_{i}\in (a;b)$ is compact and any section $\Pi _{t}\cap
\M$ is not empty. If $(a;b)$ coincide with the whole $\R{1}$ we
call $\M$ to be {\it tubular in the large}.
\end{definition}

The simplest example of a tubular in the large minimal surface in
$\R{3}$ is the catenoid. Numerous examples of two-dimensional immersed
tubular minimal surfaces in $\R{n}$ can be constructed by
the representation for such surfaces given in \cite{N}. It follows
from the papers \cite{MV}, \cite{MT2}, \cite{K} that all
$p$-dimensional tubular minimal surfaces of arbitrary codimension
and $p\geq 3$ are bounded in the $e$-direction, i.e. they can not be
tubular in the large.

To simplify the arguments we shall assume that $a=-b$. The general
case can be reduced to this one with the suitable shift of $\M$
along the $e$-axis. We set $h(m)=|\scal{x(m)}{e}| $ to be
exhausting function on $M$ provided that $\M=(M;x)$ be a minimal
surface. Moreover, $\Delta h(m)=0$ everywhere in
$M\setminus\Pi_0$ and hence, $h(m)$ is a harmonic exhausting.

\begin{theorem} %%%%%%%% 4.1
\label{th41}
Let $\M$ be a two-dimensional tubular with respect to $e$
minimal surface in $\R{n}$. Let $f$ be a harmonic function on
$M$ with the Dirichlet integral $J(t)$. If
\begin{equation}
\lim_{t \to b-0} \frac{1}{J(t)} \int_{0}^{t} \omega (y) dy = 0,
\label{new}
\end{equation}
then
$$
\liminf_{t \to b-0} \frac{\ln J(t)}{t} \geq \frac{8\pi}{S(h)}.
$$
\end{theorem}

\begin{proof} We rewrite (\ref{pro}) with applying of Proposition
\ref{prop42} as following
$$
J'(t) \geq \frac{8\pi}{S(h)} [ J(t) - Q(t) ],
$$
or,
$$
\frac{d}{dt}\left [ J(t) \exp (-8\pi t/S(h))\right ] \geq -\frac
{8\pi Q(t)}{S(h)} \exp(-8\pi t/S(h)).
$$
After integrating, we obtain
$$
J(t) \exp(-8\pi t/S(h)) \geq - \int_{0}^{t} Q(y) \exp(-8\pi y/S(h)) dy
$$
Now (\ref{new}) yields the needed assertion.
\end{proof}

\section{Some applications to the theory of minimal and $\p$-minimal
surfaces}
\label{sec5}
%%%%%%%%%%%%%%%%%%%%%%%% 5th Section %%%%%%%%%%%%%%%%%%%%%%%%%%%%%%%%%%%%%

\subsection{}
Let $M$ be an orientiable noncompact $p$-dimensional
manifold and $x : M \rightarrow \R{n}$ be a smooth isometric
immersion. Let $\nu$ be the unit normal vector field to
$\M=(M,x)$ and $A^{\nu}$ be the corresponding Weingarten map.
\par
The quantity $k_{\nu}(E) = \scal{A^{\nu}E}{E}$ is called
the curvature in $E \in T_mM$ direction.

We have already used the following characteristic property of
the minimal surfaces. Namely, they are only surfaces in $\R{n}$
with harmonic coordinate functions. On the other hand, let
$\Delta_{\p}$ be the $\p$-Laplacian of $M$ associated with the
induced metric. We fix a vector $e$ in $\R{n+1}$ and denote by
$\tau = e^{\top}/\mm{e^{\top}}$ the tangent to $\M$ direction of
the largest increasing of the coordinate function
$f_{e}(m)=\scal{e}{x(m)}$.

\begin{definition}[\cite{TVG}]
We say that a surface $\M$ is $\p$-{\it minimal} if its mean
curvature $h(m)$ and the curvature
$k_{\nu}(\tau)$ satisfy the following linear relation
\begin{equation}
h(m) = -(\p -2) k_{\nu}(\tau),
\label{kriv}
\end{equation}
at every point $m $ such that $e^{\top}(m) \neq 0$.
\end{definition}

The definition above is motivated by the following useful
analytic interpretation of condition (\ref{kriv}).

\begin{proposition}
\label{prop51}
\it
A hypersurface $\M$ is $\p$-minimal if and only if the coordinate
function $f(m)$ is $\p$-harmonic in the inner metric of $\M$.
\rm
\end{proposition}

\noindent
{\bf Proof.} We have $\gd{f} = e^{\top}$ for the gradient of
$f(m)\equiv f_e(m)$ and for a tangent vector $X$ we
obtain
$$
\nabla_X\mm{\gd{f}} = \nabla_X{\mm{\et}} = \frac{\scal{\grnd_X
\et}{\et}}{\mm{\et}}
= \frac{\scal{- \grnd_X \en}{\et}}{\mm{\et}} =
\frac{\scal{A^{\en}X}{\et}}{\mm{\et}}.
$$
Since the Weingarten map $A^{\en}$ is self-adjoint, we obtain
$$
\gd{\mm{\gd{f}}} = A^{\en}(\tau)
$$
everywhere in $\{m:\;|\gd{f}| \equiv \mm{\et} \neq 0\}$. Thus,
$$
\Delta_{\p} f = {\rm div}\mm{\gd{f}}^{\p -2}\gd{f} =
\mm{\gd{f}}^{\p -2} \Delta f + (\p -2)\scal{\nabla\mm{\gd{f}}}{\gd{f}}
\mm{\gd{f}}^{\p -3} =
$$
$$
=\mm{\gd{f}}^{\p -4} \left [ \mm{\et}^2 \Delta f + (\p -2)
\scal{A^{\en}(\et)}{\et} \right ],
$$
and applying the known \cite{KN} connection between the
Laplacian and the mean curvature
$$
%%\begin{equation}
\Delta_{\p} f = \mm{\et}^{\p -2} \scal{H +
(\p -2)B(\tau,\tau)}{\en}.
%%\label{main}
%%\end{equation}
$$
If $\M$ is a hypersurface then $\en = \scal{e}{\nu} \nu$
and the definition of $\p$-harmonicity of $f(m)$ can be
written as
$$
\bigl(h(m)+(\p-2)k_{\nu}(\tau ) \bigr) \mm{\et}^{\p -2}=0,
$$
which proves the Proposition~\ref{prop51}

%%%%%%%%%%%%%%%%%%%%%%%%% Vstavka %%%%%%%%%%%%%%%%%%%%%%%%%%%

\medskip
5.2.
We remind that for two smooth manifolds $M_1$ and $M_2$ a
mapping $F:M_1 \to M_2$ is called the mapping with {\it bounded
distortion}, or {\it quasiregular mapping}, if the Jacobian $\det
d_xF$ doesn't change the sign on $M_1$ and there exists a constant
$K\geq 1$ such that

\begin{equation}
\max |d_xF(E)|\leq K \min |d_xF(E)|,\quad x\in M,
\label{resh}
\end{equation}
where the minimum and maximum are taken over all unit tangent vectors
$E\in T_xM_1$. The smallest constant $K\equiv K(F)$ is called
the {\it distortion coefficient} of $F(x)$ (\cite{Ahlf},
\cite{Resh}).

Another well-known fact (see \cite{Os}) is that the Gauss map
$$
\gamma: M \rightarrow S^2
$$
of two-dimensional minimal surface into the standard sphere
is conformal. The following assertion was announced in \cite{TVG}
and extends that property on the class of $\p$-minimal surfaces.

\begin{theorem} %%%%%%%%% 5.1.
\label{th51}
{\it
Let $\p>1$. Then the Gauss map $\gamma$ of a two-dimensional
$\p$-minimal surface is a mapping with bounded distortion.
Moreover,
$$
K(\gamma)\leq \max (\p-1;\frac{1}{\p-1}).
$$
}
\end{theorem}

\noindent
{\bf Proof}. We observe that the tangent spaces $T_mM$ to $M$
and $T_{\gamma (m)}S^2$ to sphere $S^2$ can regarded as
canonically isomorphic ones. Really, we identify the image of the
Weingarten map $A(E)$ with $d\gamma_m(E)$.

Let us consider an arbitrary point $m$ such that $e^\top(m)\ne0$
and choose the orthonormal basis $E_1, E_2$ of the tangent space
at $m$ which one diagonalizes the symmetric Weingarten map $A^\nu$:
$$
A^\nu(E_i)=\lambda_i E_i
$$
where $\lambda_1,\,\lambda_2$ are the principal curvatures of
$\M$ at $m$. Set $\tau=e^\top/|e^\top|$. Then for some angle
$\psi \in [0;\, 2\pi )$ we have
$$
\tau = E_1 \cos\psi + E_2\sin\psi ,
$$
and by virtue of (\ref{kriv}),
$$
\scal{A^\nu (\tau)}{\tau}=\lambda_1\cos^2\psi + \lambda_2\sin^2\psi
= - \frac{1}{\p-2}(\lambda_1+\lambda_2).
$$
Thus,
$$
\lambda_1 = -\lambda_2\;\frac{1+(\p-2)\sin^2\psi}{1+(\p-2)\cos^2\psi}.
$$

We observe that the last identity yields negativity of Jacobian
$\det(d_m\gamma )= \lambda_1\lambda_2$ at $m$. Now, using the
extremal properties of the quadratic forms we obtain the value
of the distortion coefficient $K$ at $m$,
$$
K_m=\max_{\psi}\{q ; \frac{1}{q}\} ,
\qquad q=\frac{1+(\p-2)\sin^2\psi}{1+(\p-2)\cos^2\psi}.
$$
An easy computation of the maximal value of the right part over
all admissible angles $\psi$ yields the required inequality.

If $e^\top(m)=0$, then one can show (see \cite{TVG})
that Weingarten map $A^\nu$ is the identical null and
(\ref{resh}) is trivial.

\begin{corollary}
\label{cor51}
{\it
For $\p>1$ two-dimensional planes are only entire
$\p$-minimal graphs in $\R{3}$.
}
\end{corollary}

\noindent
{\bf Proof.}
This version of the well-known S.N.Bernstein's result is a direct
consequence of the previous theorem and the theorem of L.Simon
\cite{SL} on the two-dimensional entire graphs with
quasiconformal Gauss map.
%%%%%%%%%%%%%%%%%% The end of Vstavka

\subsection{} Let $\p = 2$, $M$ be a two-dimensional manifold and $f(m)$ be a
subharmonic function on $M$ having $N$ different regular
asymptotic tracts $\D_1, \D_2, \ldots, \D_N$. Let $h(m)$ be an
exhausting function on $M$ and
$$
\theta (t) = \; \int\limits_{\Sigma_h(t)} \mm{\gd{h}},
$$
is the flow of $h(m)$ through $\Sigma_h(t)$. Then we have from Lemma~\ref{lem25}
and inequalities (\ref{estLam}), (\ref{Nmin}) for any $t_1 < t_2 < t_3$
such as $t_1 > \max_{i\leq N} h(\D_i)$ :
\smallskip
\begin{equation}
\frac{N}{4}\min_{i\leq N} \J{2}{f}{\D_i\cap B_h(t_1)}
\leq \left [ \int\limits_{t_2}^{t_3} \frac{dt}{\mu^2(t)\theta (t)}\right
]^{-1}
\exp\left ( -2\pi N \int\limits_{t_1}^{t_2}\frac{dt}{\theta (t)}\right ),
\label{pp}
\end{equation}
where $\mu (t) = \max_{m\in \Sigma_h(t)} \mm{f(m)}$.

Now we suppose that $x : M \to \R{n}$ is a proper minimal
immersion. Then $x_k(m), 1 \leq k \leq n$, are harmonic
functions on $M$ and Lemma~\ref{lem11} yields that $h(m) = \mm{x(m)}$ is
an exhausting function. We observe that
$$
\gd{h(m)} = (\grnd{\mm{x(m)}})^{\top} = \frac{x^{\top}(m)}{\mm{x(m)}},
$$
and hence,
$$
\mm{\gd{h(m)}} = \frac{\mm{x^{\top}(m)}}{\mm{x(m)}} \leq 1
$$
everywhere in $M$, and after applying Cauchy's inequality,
$$
\ln {b\over a} \leq \left ( \int\limits_{a}^{b} {\theta (t)
dt\over t^2} \right )^{1/2}\left ( \int\limits_{a}^{b}
{dt\over \theta (t)} \right )^{1/2},
$$
we obtain

\begin{equation}
\int\limits_{a}^{b} {dt\over \theta (t)} \geq (\ln \frac{b}{a})^2
\left [ \int\limits_{a}^{b}\frac{dt}{t^2} \; \int\limits_{\Sigma_t}
\mm{\gd{h}} \right ]^{-1} =
(\ln \frac{b}{a})^2
\left [ \int\limits_{B_h(b)\setminus B_h(a)}\frac{1}{\mm{x}^2} \right
]^{-1}
\label{nerav}
\end{equation}

Substituting in (\ref{pp}) the coordinate function $x_k(m)$ instead
of $f(m)$ arrive at
\begin{equation}
\min_{1\le i\le N} \int\limits_{B_h(t_1)\cap \D_i} \mm{\gd{x_k}}^2 \leq
\frac{4\, t_3^2\, (V(t_3) - V(t_2))}{N\, \ln^2 \,(t_3/t_2)}
\exp \left (-\frac{2\pi N\ln^2(t_2/t_1)}{V(t_2) - V(t_1)}\right ).
\label{ve}
\end{equation}
Here
$$
V(t) = \int_{B_h(t)}\frac{1}{|x|^2}.
$$

This quantity have been used by the authors in the paper \cite{MT1}
for estimation of the extremal length of a family of curves on
minimal surface. In fact, the asymptotic behavior of $V(t)$ at
infinity can describe in terms of the integral-geometrical
invariants of minimal submanifolds. In recent paper \cite{T1}
the following property of $V(t)$ has been established. Let $\M$
be a $p$-dimensional properly immersed minimal surface
in $\R{n}$ such that
\begin{equation}
V(t)\equiv\int_{M(t)}\frac{1}{|x(m)|^p}=O(\ln t), \quad\mbox{as
$t\to \infty$}
\label{area}
\end{equation}
where $M(t)=\{m:1<|x(m)|<t\}$. Then \cite{T1} there exist the limits
\begin{equation}
V_p(\M)\equiv\lim_{t\to\infty}\frac{V(t)}{{\omega_{p}}\ln
t}=\lim_{t\to\infty}\frac{p\;{\rm Area}_p(M(t))}{\omega_p t^p},
\label{volume}
\end{equation}
where $\omega_p$ is the $(p - 1)$-dimensional measure of a
unit sphere $S^{p-1}(1)$.

\medskip

{\bf Remark 2.} In particularly (see \cite{T1}), if $\M$ is a
two-dimensional properly immersed surface of finite total
Gaussian curvature then $V_2(\M)=\ell$, where $\ell$ is the
number of ends of $\M$. Moreover, it was shown in \cite{T1}
that properly immersed $p$-dimensional minimal submanifolds with
$V_p(\M)<\infty$ have the finite number of topological ends. On
the other hand, for the two-dimensional helicoid one holds
$V_2(\M)=+\infty$ while its Eulerian characteristic and the
number $\ell$ are finite. We don't also know whether
$V_2(\M)<+\infty$ yields finiteness of the total Gaussian
curvature of a minimal surface $\M$.

In this paragraph we study relations between $V_2(\M)$ and the
number of humps which can be cut off from a two-dimensional
minimal surface by a system of hyperplanes. In this connection,
we notice that even if $\M$ be the two-dimensional
catenoid then there are planes which don't cut off any humps from the
surface. Really, it is sufficient to consider a plane which is
orthogonal to the axis of the catenoid. Therefore, we need the
following

\medskip
{\bf Definition 12.} A direction $e\in\R{n}$ is called
{\it regular}\/ for a surface $\M$ if $\M$ doesn't contain in a
hyperplane orthogonal to $e$ and all sections of $\M$ by
hyperplanes $\Pi\bot e$ don't contain compact components (here
and henceforth we mean by "sections of $\M$ etc" the words "the
preimage of sections of $x(M)$ etc" if the surface $\M$ doesn't
embedded).

For such two hyperplanes $\Pi_1$ and $\Pi_2$ we denote by
$N(\Pi_1;\Pi_2)$ the number of components of $\M$ lying outward
of the slab with a boundary $\Pi_1\cup \Pi_2$.

We notice, that a hump must has noncompact boundary.

\begin{theorem} %%%%%%% 5.2.
\label{th52}
Let $\M$ be a two-dimensional properly immersed minimal surface in
$\R{n}$ with finite projective $2$-volume $V_2(\M)$.
Let $e$ be a regular direction. Then for any hyperplanes $\Pi_1$
and $\Pi_2$ orthogonal to $e$
\begin{equation}
N(\Pi_1;\Pi_2) \; \leq \; 2 V_2(\M).
\label{enqu}
\end{equation}
\end{theorem}

\noindent
{\bf Proof.}
Really, let $x_1$ be corresponding to $e$ coordinate function
$x_1(m)=\scal{x(m)}{e}$. Then $\M$ doesn't contain in
any hyperplane which is orthogonal to $e$ and it follows that
$x_1(m)$ does not constant. Without loss of generality we can arrange that $\Pi_i$
is defined by $x_1=(-1)^ia$ for some $a>0$ and consider the
subharmonic function $f(m)=|x_1(m)|$. We fix $t_1>a$ and denote by
$J>0$ the left part of (\ref{ve}). Given
$\varepsilon > 0$, we find $t_2$ to be sufficient large such that
$$
V(t) < 2\pi (\q + \varepsilon) \ln t.
$$
Hence, for any $t > t_2$ and from (\ref{pp}),
$$
J < \frac{4 t_3^2 \, t_2^{-N/(\q +\varepsilon)}}{N \ln
(t_3/t_2)},
$$
where $t_3>t_2$ and $N=N(\Pi_1;\Pi_2)$.

We choose now $t_3 = 2 t_2$. Then $t_2 \to\infty $ gives
$$
2 - \frac{N}{\q + \varepsilon} \; \geq 0,
$$
and (\ref{enqu}) is proved.

\medskip
5.4. In this paragraph we deal with the minimal surfaces $\M$ of
finite topological type. This means that $\M$ is realized by a
minimal immersion of some compact manifold $M$ of genus $g$ with
a finite number $\ell$ of points removed. The last points is called
the {\it ends} of $M$.

{\bf Definition 13.} Let $f(m)$ be a harmonic function on $M$. Due
to \cite{Mr} we define the {\it index} of the function $f(m)$ at a
critical point $m_0\subset {\mathcal Z}(f)$ to be the number
$$
{\rm ind}(m_0) = \frac{\sigma}{2} - 1,
$$
where $\sigma$ is the number of splitting continua of the set
$\{m\in M: f(m) = f(m_0), m\ne m_0 \mbox{and $m$ is sufficiently
near to $m_0$} \}$.

It has been shown in \cite{Mr} that ${\rm ind}(m_0)$ is a positive
integer provided that $f$ is not a constant.

\begin{theorem} %%%%%%%%%%% 5.3.
\label{th53}
Let $\M $ be a two-dimensional properly immersed minimal surface
in $\R{n}$ of finite topological type, $e$ is a regular
direction and $x_1$ be the corresponding coordinate function. If
$V_2(\M)<+\infty$ then the number of critical points $\{ a_i \}$ of
$x_1(m)$ is finite. Moreover,
\begin{equation}
\sum_{j} {\rm ind}(a_j) \leq \q - \chi (M),
\label{pom}
\end{equation}
where $\chi (M)$ is the Eulerian characteristic of $M$.
\end{theorem}

{\bf Remark 3.} When $M$ is homeomorphic to a sphere with $\nu$
points removed this assertions was proved in \cite{M2}.

\medskip
\noindent {\bf Proof.} We denote by ${\mathcal N} (t)$ the number
of components of $\{ m \in M : x_1(m) > t \}$. By virtue of the
maximum principle for the coordinate functions of minimal
surfaces, ${\mathcal N} (t)$ is nondecreasing integer-valued
function for $t \in \R{}$, and for any critical value $a_i$ we
have
\begin{equation}
\lim_{t\to x_1(a_i)+0} {\mathcal N} (t) - \lim_{t\to x_1(a_i)-0}
{\mathcal N} (t) \geq 1. \label{pop}
\end{equation}

Really, the level set $\{ m \in M :\; x_1(m) = a_i \}$ is a
union of continua $\gamma$ and by virtue of assumption of
the regularity of the coordinate function, none of $\gamma$'s is
compact. We write
$$
{\mathcal N}(+\infty)\equiv \lim_{t\to +\infty} {\mathcal N}(t).
$$

We also observe that finiteness of the projective volume of $\M $
yields parabolic conformal type of $\M $ (see \cite{MT1}). Thus,
by the Phragmen-Lindel\"of theorem applied to the harmonic
function $x_1(m)$ we can suppose that all humps of $\{m\in
M|x_1(m)>t\}$ are unbounded at the positive $x_1$-direction. Then
finiteness of the number of critical points $a_i$ of $x_1$ follows
from (\ref{pop}) and finiteness of the number ${\mathcal N}$ of
asymptotic tracts of $|x_1(m)|$. The last property is a
consequence of (\ref{enqu}).

Our proof of (\ref{pom}) is an appropriate generalization of the
corresponding estimate in \cite{M2} for zero-genus minimal
surfaces.

We choose $c$ to be a positive number which is greater than the
absolute value of any critical value $x_1(a_j)$.
It follows from conformality of the Gauss map of a twodimensional
minimal surface that it has locally finite multiplicity. This
property provides (see \cite[Theorem 18.5.4]{Bak}) existence of
the {\it regular} part $M'$ of $M$, i.e. such a smooth compact
submanifold of $M$ (with nonempty boundary) which satisfies the
following conditions

\begin{enumerate}
\item
the topological type of $M'$ coincides with the same of $M$:
$\chi(M')=\chi(M)$;
\item
every boundary component of $M$ is a finite family of
alternating $x_1(m)$-level-curves $\gamma_i$ ($x_1=\pm c$) and
the gradient-curves $\Gamma_i$ (i.e. the curves of the most
increasing of $x_1(m)$);
\item
the number $\ell$ of all components of $\partial M'$ coincides
with the number of ends of $M$.
\end{enumerate}

From regularity of the $e$-direction follows that each
$\gamma_i$ is a simple open curve. As a consequence, every
component of $\partial M'$ contains of an even number of
$\gamma_i$ and the same number of $\Gamma_j$.

We consider the decomposition of the following two parts of $M$
into open components
$$
\{ x_1(m) > c \} = \cup_{i=1}^{N^{+}} \OO_i^{+} \quad,
\quad \{ x_1(m) < - c \} = \cup_{j=1}^{N^{-}} \OO_j^{-}.
$$

Since the set $\{|x_1(m)|>c\}$ contains no critical
points, each $\OO_i^\pm$ is homotopically equivalent to
a two-dimensional disk. Furthermore, one can determinate a
bijection between the arcs $\gamma_i$ and the
components $\OO_i^{\pm}$. It follows from (\ref{enqu}) and the
alternating property of $\gamma_i$ and $\Gamma_j$ that
$$
N^+ = N^- = {\mathcal N}(+\infty)=\frac{{\mathcal N}}{2}\leq \q.
$$

Let $M_0 = M'\# M'$ be the result of pasting together of two
copies of $M'$ along $\partial M'$ and consequent contracting of
each curve $\gamma_i$ and its copy into a point $G_i$.
Then (see \cite{Sp}, Section~5, ex.~5), the genus $g_0$ of $M_0$
can be expressed by
\begin{equation}
g_0 = 2g + \ell - 1,
\label{gzero}
\end{equation}
where $\ell$ is a number of the components of $\partial M'$.

By virtue of the definition of $M_0$ the coordinate function
$x_1(m)$ can be canonically lifted up to the well-defined on
$M_0$ function $f$. Then $f(\xi)$ is smooth everywhere on
$M_0\setminus G$, where $G=\{G_i\}$. Moreover, all of $G_j$'s
are the points of strong maximum and minimum of $f(\xi)$.

For the Eulerian characteristic we have from (\ref{gzero})
\begin{equation}
\chi (M_0)=2-2g_0=2(2-2g-\ell) =2\chi (M')=2 \chi (M).
\label{pos}
\end{equation}

We can lift by a natural way the gradients field $\nabla
x_1(m)$ up to the continuous vector field $X(m)\equiv \nabla f$
on $M_0\setminus G$. Then the set of
singular points of $X(m)$ consists of: a) the critical
points $\{a_i\}^k_{i=1}$ and its doubles $\{ a_i^{\star} \}^k_{i=1}$;
and b) the singular points $G_j$.

Let ${\rm ind}_{\xi}X$ denotes the rotational number ({\it
index}) of the field $X$ at $\xi$ \cite[\S 5]{Mil}. Then at the
extreme points we have
\begin{equation}
{\rm ind}_{G_j}X=1.
\label{inde}
\end{equation}

On the other hand, by virtue of the definition of the index of
harmonic function (see also \cite[\S6]{Mil}) it is easily seen
\begin{equation}
{\rm ind}_{a_i} X= {\rm ind}_{a_i^\star} X\equiv {\rm ind}_{a_i}
(\nabla x_1) =- {\rm ind}(a_i).
\label{inde2}
\end{equation}
Indeed, the last identity follows from analysis of a
sufficient small neighbourhood of the critical point $a_i$ of a
harmonic function.

Thus, we can apply the Poincare-Hopf theorem and from
(\ref{inde}), (\ref{inde2}) we obtain
$$
\chi(M_0)=\sum {\rm ind}_{a_i}X+\sum {\rm ind}_{a_i^\star}X
+\sum_{j=1}^{\N}{\rm ind}_{G_j}X=
$$
$$
-2\;\sum {\rm ind}(a_i)+\N.
$$
Hence, by (\ref{pos}) we arrive at
$$
\sum {\rm ind}(a_i)=\frac{\N}{2}-\chi(M)\leq\q-\chi(M),
$$
and the theorem is proved.

\medskip
5.5. To present some applications of these theorems we first
suppose that $\M$ is a properly immersed plane in $\R{3}$, i.e.
$\chi (M) = \chi (\R{2}) = 1$.

{\bf Definition 14.} Let $V$ be a 2-dimensional plane in
$\R{3}$ and $\pi_V:\R{3}\to V$ be the orthogonal projection. A
surface $\M=(M;u)$ is said to be {\it proper} with respect to
$V$ if for any sequence $\{m_k\}\subset M$ without accumulation
points in $M$ the sequence $\pi_V\circ u(m_k)$ hasn't any
accumulation point in $V$. Other words, the composition
$\pi_V\circ u$ is proper mapping.

Given a plane $V$ in $\R{3}$ denote by $n(v)$ the algebraic
multiplicity of the orthogonal projection of $\M$ onto $V$ at
$v$.

\begin{theorem} %%%%% 5.4
\label{th54}
Let $\M$ be a two-dimensional minimal properly immersed plane in
$\R{3}$ and $\M$ be proper with respect to a two-dimensional
plane $V$. Then either $\M$ is a plane, or
\begin{equation}
\liminf_{R\to\infty} \frac{1}{\ln R} \int_{1}^{R} \frac{dt}{t}
\int_{0}^{2\pi} n(t e^{i\theta}) d\theta \geq 8.
\label{pol}
\end{equation}
\end{theorem}

\noindent
{\bf Proof.} In \cite{MT1} (see also \cite{T1}) was proved that
\begin{equation}
\q \le \liminf_{R\to\infty} \frac{1}{4\ln R} \int_{1}^{R} \frac{dt}{t}
\int_{0}^{2\pi} n(t e^{i\theta}) d\theta.
\label{pot}
\end{equation}
Hence $\M$ has finite projective volume $\q$. Let $\sigma :
M \rightarrow {\rm S}^2$ be the Gauss map and $\M$ be different
from a plane. Then $\sigma (M)$ has nonempty interior ${\rm int}
\sigma (M)$. We choose arbitrary $\nu\in{\rm int}\sigma (M)$. Thus
the corresponding coordinate function
$$
f(m) = \scal{x(m)}{\nu} = \nu_1 x_1(m) + \nu_2 x_2(m) + \nu_3 x_3(m)
$$
has at least one critical point $m_0$ corresponding to
$\sigma (m_0) = \nu$. We notice, that by virtue of our
assumption $\chi(\M)=1$, $\nu$ is a regular direction provided
$\M$ doesn't consist in a hyperplane. Therefore, we have
\begin{equation}
{\rm ind}(m_0)\geq1.
\label{por}
\end{equation}

Now our statement follows from $\chi (\M )= 1$, (\ref{por}),
(\ref{pot}) and (\ref{pom}).

As a consequence of this theorem we obtain a new proof of the
Bernstein`s theorem

\begin{corollary}
\label{cor52}
The only entire graphs of minimal surfaces are the planes.
\end{corollary}

\noindent
{\bf Proof}. We observe that $n(v) \equiv 1 $ for a graph
and therefore
$$
\liminf_{R\to\infty} \frac{1}{\ln R} \int_{1}^{R} \frac{dt}{t}
\int_{0}^{2\pi} n(t e^{i\theta}) d\theta = 2 \pi < 8.
$$

Another consequence gives a quantitative form of the previous result.

\begin{corollary}
\label{cor53}
\it
Let $\M$ be an immersed minimal surface of finite topological
type which is homeomorphic to a compact Riemannian surface of
genus $g$ with $l$ points removed. Then in the given above notations:
$$
\liminf_{R\to\infty} \frac{1}{\ln R}\int_{1}^{R} \frac{dt}{t}
 \int_{0}^{2 \pi} n( t e^{i\theta})
d\theta \geq 2(l + 3 - g).
$$
\rm
\end{corollary}

%%%%%%%%%%%%%%%%%%%%%%%%%%%%%%%%%%%%%%%%%%%%%%%%%%%%%%%%%%%%%%%%%%%%%%%

\end{document}